\documentclass[11pt,leqno]{article}
\usepackage{amssymb}
\usepackage[mathscr]{eucal}
\usepackage{amsmath,amssymb,latexsym,theorem,bbm}
\usepackage{color}
\usepackage{appendix}
\usepackage{listings}
\usepackage{graphicx}

\newcommand{\NN}{\mathbb{N}}

\newcommand{\RR}{\mathbb{R}}

\newcommand{\ZZ}{\mathbb{Z}}

\newcommand{\bA}{{\boldsymbol{A}}}

\newcommand{\bB}{{\boldsymbol{B}}}

\newcommand{\bC}{{\boldsymbol{C}}}

\newcommand{\bD}{{\boldsymbol{D}}}

\newcommand{\bI}{{\boldsymbol{I}}}

\newcommand{\bM}{{\boldsymbol{M}}}

\newcommand{\bQ}{{\boldsymbol{Q}}}
\newcommand{\bS}{{\boldsymbol{S}}}

\newcommand{\bv}{{\boldsymbol{v}}}

\newcommand{\bZ}{{\boldsymbol{Z}}}

\newcommand{\bfeta}{{\boldsymbol{\eta}}}

\newcommand{\bSigma}{{\boldsymbol{\Sigma}}}
\newcommand{\bXi}{{\boldsymbol{\Xi}}}
\newcommand{\bzero}{{\boldsymbol{0}}}

\newcommand{\cF}{{\mathcal F}}

\newcommand{\cL}{{\mathcal L}}

\newcommand{\cN}{{\mathcal N}}

\newcommand{\dd}{\mathrm{d}}
\newcommand{\ee}{\mathrm{e}}

\newcommand{\LSE}{\mathrm{LSE}}
\newcommand{\LSED}{\mathrm{LSE,D}}

\DeclareMathOperator*{\argmin}{arg\,min}

\newcommand{\EE}{\operatorname{\mathbb{E}}}
\newcommand{\PP}{\operatorname{\mathbb{P}}}

\newcommand{\ha}{\widehat{a}}
\newcommand{\hb}{\widehat{b}}

\newcommand{\halpha}{\widehat{\alpha}}
\newcommand{\hbeta}{\widehat{\beta}}

\newcommand{\tW}{\widetilde{W}}

\newcommand{\oE}{\overline{E}}

\renewcommand{\mid}{\,|\,}

\renewcommand{\leq}{\leqslant}
\renewcommand{\geq}{\geqslant}

\newcommand{\stoch}{\stackrel{\PP}{\longrightarrow}}
\newcommand{\distr}{\stackrel{\cL}{\longrightarrow}}
\newcommand{\distre}{\stackrel{\cL}{=}}

\newcommand{\as}{\stackrel{{\mathrm{a.s.}}}{\longrightarrow}}

\newcommand{\nT}{{\lfloor nT\rfloor}}

\newcommand{\proofend}{\hfill\mbox{$\Box$}}

\numberwithin{equation}{section}

\theoremstyle{change} \theorembodyfont{\em}
\newtheorem{Lem}{Lemma.}[section]
\newtheorem{Thm}[Lem]{Theorem.}
\newtheorem{Pro}[Lem]{Proposition.}
\newtheorem{Cor}[Lem]{Corollary.}
\newtheorem{Def}[Lem]{Definition.}

\theorembodyfont{\rm}
\newtheorem{Rem}[Lem]{Remark.}

\begin{document}

\begin{center}
 {\bfseries\Large
  Least squares estimation for the subcritical Heston model \\ based on continuous time observations} \\[5mm]
 {\sc\large
  M\'aty\'as $\text{Barczy}^{*,\diamond}$,
  Bal\'azs $\text{Nyul}^{**}$ \ and
  \ Gyula $\text{Pap}^{***}$}
\end{center}

\vskip0.2cm

\noindent
 * MTA-SZTE Analysis and Stochastics Research Group,
   Bolyai Institute, University of Szeged,
   Aradi v\'ertan\'uk tere 1, H--6720 Szeged, Hungary.

\noindent
 ** Faculty of Informatics, University of Debrecen,
    Pf.~12, H--4010 Debrecen, Hungary.

\noindent
 *** Bolyai Institute, University of Szeged,
      Aradi v\'ertan\'uk tere 1, H--6720 Szeged, Hungary.

\noindent e--mails: barczy@math.u-szeged.hu (M. Barczy),
                   nyul.balazs@inf.unideb.hu (B. Nyul),\\
                   papgy@math.u-szeged.hu (G. Pap).

\noindent $\diamond$ Corresponding author.

\renewcommand{\thefootnote}{}
\vspace*{0.2cm}
\footnote{M\'aty\'as Barczy is supported by the J\'anos Bolyai Research Scholarship of the Hungarian Academy
 of Sciences.}

\vspace*{-5mm}


\begin{abstract}
We prove strong consistency and asymptotic normality of least squares estimators for the subcritical Heston model
 based on continuous time observations.
We also present some numerical illustrations of our results.
\end{abstract}

\section{Introduction}

Stochastic processes given by solutions to stochastic differential equations (SDEs) have been frequently applied in financial
 mathematics.
So the theory and practice of stochastic analysis and statistical inference for such processes are important topics.
In this note we consider such a model, namely the Heston model
 \begin{align}\label{Heston_SDE}
  \begin{cases}
   \dd Y_t = (a - b Y_t) \, \dd t + \sigma_1 \sqrt{Y_t} \, \dd W_t , \\
   \dd X_t = (\alpha - \beta Y_t) \, \dd t
             + \sigma_2  \sqrt{Y_t}
               \bigl(\varrho \, \dd W_t
                     + \sqrt{1 - \varrho^2} \, \dd B_t\bigr) ,
  \end{cases} \qquad t \geq 0 ,
 \end{align}
 where \ $a > 0$, \ $b, \alpha, \beta \in \RR$, \ $\sigma_1 > 0$,
 \ $\sigma_2 > 0$, \ $\varrho \in (-1, 1)$, \ and \ $(W_t, B_t)_{t\geq0}$ \ is a
 2-dimensional standard Wiener process, see Heston \cite{Hes}.
For interpretation of \ $Y$ \ and \ $X$ \ in financial mathematics, see, e.g., Hurn et al. \cite[Section 4]{HurLinMcC},
 here we only note that \ $X_t$ \ is the logarithm of the asset price at time \ $t$ \
 and \ $Y_t$ \ its volatility for each \ $t\geq 0$.
\ The first coordinate process \ $Y$ \ is called a Cox-Ingersoll-Ross (CIR) process
 (see Cox, Ingersoll and Ross \cite{CoxIngRos}), square root process or Feller process.

Parameter estimation for the Heston model \eqref{Heston_SDE} has a long history, for a short survey
 of the most recent results, see, e.g., the introduction of Barczy and Pap \cite{BarPap}.
The importance of the joint estimation of \ $(a, b, \alpha, \beta)$ \ and not only of \ $(a,b)$ \ stems from the fact that
 \ $X_t$ \ is the logarithm of the asset price at time \ $t$ \ having high importance in finance.
In fact, in Barczy and Pap \cite{BarPap}, we investigated asymptotic properties of maximum likelihood estimator of \ $(a, b, \alpha, \beta)$ \ based on continuous time observations \ $(X_t)_{t\in[0,T]}$, $T>0$.
\ In Barczy et al.~\cite{BarPapSza} we studied asymptotic behaviour of conditional least squares
 estimator of \ $(a, b, \alpha, \beta)$ \ based on discrete time observations
 \ $(Y_i,X_i)$, $i=1,\ldots,n$, \ starting the process from some known non-random initial value
 \ $(y_0,x_0)\in(0,\infty)\times \RR$.
\ In this note we study least squares estimator (LSE) of \ $(a, b, \alpha, \beta)$ \ based
 on continuous time observations \ $(X_t)_{t\in[0,T]}$, \ $T > 0$,
 \ starting the process \ $(Y,X)$ \ from some known initial value
 \ $(Y_0,X_0)$ \ satisfying \ $\PP(Y_0\in(0,\infty))=1$.
\ The investigation of the LSE of \ $(a,b,\alpha,\beta)$ \ based on continuous time observations \ $(X_t)_{t\in[0,T]}$, $T>0$, \
 is motivated by the fact that the LSEs of \ $(a,b,\alpha,\beta)$ \ based on appropriate discrete time observations converge
 in probability to the LSE of \ $(a,b,\alpha,\beta)$ \  based on continuous time observations \ $(X_t)_{t\in[0,T]}$, $T>0$, \
 see Proposition \ref{LEMMA_LSE_exist}.
We do not suppose that the process \ $(Y_t)_{t\in[0,T]}$ \ is observed, since it
 can be determined using the observations \ $(X_t)_{t\in[0,T]}$ \ and the initial value
 \ $Y_0$, \ which follows by a slight modification of Remark 2.5 in Barczy and Pap \cite{BarPap}
 (replacing \ $y_0$ \ by \ $Y_0$).
\ We do not estimate the parameters \ $\sigma_1$, \ $\sigma_2$ \ and \ $\varrho$,
 \ since these parameters could ---in principle, at least--- be determined
 (rather than estimated) using the observations \ $(X_t)_{t\in[0,T]}$ \
 and the initial value \ $Y_0$, \ see Barczy and Pap \cite[Remark 2.6]{BarPap}.
We investigate only the so-called subcritical case, i.e., when \ $b>0$,
 \ see Definition \ref{Def_criticality}.

In Section \ref{Prel} we recall some properties of the Heston model \eqref{Heston_SDE}
 such as the existence and uniqueness of a strong solution of the SDE \eqref{Heston_SDE},
 the form of conditional expectation of \ $(Y_t,X_t)$, \ $t \geq 0$, \ given the past of the process up to time \ $s$ \ with \ $s \in[0, t]$, \ a classification
 of the Heston model and the existence of a unique stationary distribution and ergodicity
 for the first coordinate process of the SDE \eqref{Heston_SDE}.
Section \ref{section_LSE_continuous} is devoted to derive a LSE of \ $(a, b, \alpha, \beta)$ \ based
 on continuous time observations \ $(X_t)_{t\in[0,T]}$, \ $T > 0$, \ see Proposition \ref{LEMMA_LSE_exist}.
We note that Overbeck and Ryd\'en \cite[Theorems 3.5 and 3.6]{OveRyd} have already proved the
 strong consistency and asymptotic normality of the LSE of \ $(a,b)$ \ based on continuous time observations
 \ $(Y_t)_{t\in[0,T]}$, $T>0$, \ in case of a subcritical CIR process \ $Y$ \ with an initial value
 having distribution as the unique stationary distribution of the model.
Overbeck and Ryd\'en \cite[page 433]{OveRyd} also noted that (without providing a proof) their results are valid for an
 arbitrary initial distribution using some coupling argument.
In Section \ref{section_cons_asy} we prove strong consistency and asymptotic normality of the LSE
 of \ $(a, b, \alpha, \beta)$ \ introduced in Section \ref{section_LSE_continuous},  so our results
 for the Heston model \eqref{Heston_SDE} in Section \ref{section_LSE_continuous} can be considered as generalizations
 of the corresponding ones in Overbeck and Ryd\'en \cite[Theorems 3.5 and 3.6]{OveRyd}
 with the advantage that our proof is presented for an arbitrary initial value \ $(Y_0,X_0)$ \
 satisfying \ $\PP(Y_0\in(0,\infty))=1$, \ without using any coupling argument.
The covariance matrix of the limit normal distribution in question depends on the unknown parameters
  \ $a$ \ and \ $b$ \ as well, but somewhat surprisingly not on \ $\alpha$ \ and \ $\beta$.
\ We point out that our proof of technique for deriving the asymptotic normality of the LSE in question
 is completely different from that of Overbeck and Ryd\'en \cite{OveRyd}.
We use a limit theorem for continuous martingales (see, Theorem \ref{THM_Zanten}),
 while Overbeck and Ryd\'en \cite{OveRyd} use a limit theorem for ergodic processes due to Jacod and Shiryaev
 \cite[Theorem VIII.3.79]{JSh} and the so-called Delta method (see, e.g., Theorem 11.2.14 in Lehmann and Romano \cite{LehRom}).
We also remark that the approximation in probability of the LSE of \ $(a,b,\alpha,\beta)$ \ based on continuous time observations
 \ $(X_t)_{t\in[0,T]}$, $T>0$, \ given in Proposition \ref{LEMMA_LSE_exist} is not at all used for proving the asymptotic behaviour
 of the LSE in question as \ $T\to\infty$ \ in Theorems \ref{Thm_LSE_cons} and \ref{Thm_LSE}.
Further, we mention that the covariance matrix of the limit normal distribution in Theorem 3.6 in Overbeck and Ryd\'en
 \cite{OveRyd} is somewhat complicated, while, as a special case of our Theorem \ref{Thm_LSE},
 it turns out that it can be written in a much simpler form by making a simple reparametrization
 of the SDE (1) in Overbeck and Ryd\'en \cite{OveRyd}, estimating \ $-b$ \ instead of \ $b$ \ (with the notations
 of Overbeck and Ryd\'en \cite{OveRyd}), i.e., considering the SDE \eqref{Heston_SDE} and estimating \ $b$ \ (with our notations),
 see Corollary \ref{Rem_Ove_Ryd}.
Section \ref{section_numerical} is devoted to present some numerical illustrations of our results
 in Section \ref{section_cons_asy}.

\section{Preliminaires}
\label{Prel}

Let \ $\NN$, \ $\ZZ_+$, \ $\RR$, \ $\RR_+$, \ $\RR_{++}$, \ $\RR_-$ \ and
 \ $\RR_{--}$ \ denote the sets of positive integers, non-negative integers,
 real numbers, non-negative real numbers, positive real numbers, non-positive
 real numbers and negative real numbers, respectively.
For \ $x , y \in \RR$, \ we will use the notation \ $x \land y := \min(x, y)$.
\ By \ $\|x\|$ \ and \ $\|A\|$, \ we denote the Euclidean norm of a vector
 \ $x \in \RR^d$ \ and the induced matrix norm of a matrix
 \ $A \in \RR^{d \times d}$, \ respectively.
By \ $\bI_d \in \RR^{d \times d}$, \ we denote the $d$-dimensional unit matrix.

Let \ $\bigl(\Omega, \cF, \PP\bigr)$ \ be a probability space
 equipped with the augmented filtration \ $(\cF_t)_{t\in\RR_+}$ \ corresponding to
 \ $(W_t,B_t)_{t\in\RR_+}$ \ and a given initial value \ $(\eta_0,\zeta_0)$ \ being independent
 of \ $(W_t,B_t)_{t\in\RR_+}$ \ such that \ $\PP(\eta_0\in\RR_+)=1$, \
 constructed as in Karatzas and Shreve \cite[Section 5.2]{KarShr}.
Note that \ $(\cF_t)_{t\in\RR_+}$ \ satisfies the usual conditions, i.e., the
 filtration \ $(\cF_t)_{t\in\RR_+}$ \ is right-continuous and \ $\cF_0$
 \ contains all the $\PP$-null sets in \ $\cF$.

By \ $C^2_c(\RR_+\times\RR, \RR)$ \ and \ $C^{\infty}_c(\RR_+\times\RR, \RR)$,
 \ we denote the set of twice continuously differentiable real-valued
 functions on \ $\RR_+\times\RR$ \ with compact support, and the set of
 infinitely differentiable real-valued functions on \ $\RR_+\times\RR$ \ with
 compact support, respectively.

The next proposition is about the existence and uniqueness of a strong
 solution of the SDE \eqref{Heston_SDE}, see, e.g., Barczy and Pap
 \cite[Proposition 2.1]{BarPap}.

\begin{Pro}\label{Pro_Heston}
Let \ $(\eta_0, \zeta_0)$ \ be a random vector independent of
 \ $(W_t, B_t)_{t\in\RR_+}$ \ satisfying \ $\PP(\eta_0 \in \RR_+) = 1$.
\ Then for all \ $a \in \RR_{++}$, \ $b, \alpha, \beta \in \RR$,
 \ $\sigma_1, \sigma_2 \in \RR_{++}$, \ and \ $\varrho \in (-1, 1)$, \ there is a
 pathwise unique strong solution \ $(Y_t, X_t)_{t\in\RR_+}$ \ of the SDE
 \eqref{Heston_SDE} such that \ $\PP((Y_0, X_0) = (\eta_0, \zeta_0)) = 1$ \ and
 \ $\PP(\text{$Y_t \in \RR_+$ \ for all \ $t \in \RR_+$}) = 1$.
\ Further, for all \ $s, t \in \RR_+$ \ with \ $s \leq t$,
 \begin{align}\label{Solutions}
  \begin{cases}
   Y_t = \ee^{-b(t-s)} Y_s
                + a \int_s^t \ee^{-b(t-u)} \, \dd u
                + \sigma_1
                  \int_s^t \ee^{-b(t-u)} \sqrt{Y_u} \, \dd W_u , \\
   X_t = X_s + \int_s^t (\alpha - \beta Y_u) \, \dd u
         + \sigma_2
           \int_s^t
            \sqrt{Y_u} \, \dd (\varrho W_u + \sqrt{1 - \varrho^2} B_u) .
  \end{cases}
 \end{align}
\end{Pro}

Next we present a result about the first moment and the conditional moment of \ $(Y_t, X_t)_{t\in\RR_+}$, \
 see Barczy et al. \cite[Proposition 2.2]{BarPapSza}.

\begin{Pro}\label{Pro_moments}
Let \ $(Y_t, X_t)_{t\in\RR_+}$ \ be the unique strong solution of the SDE
 \eqref{Heston_SDE} satisfying \ $\PP(Y_0 \in \RR_+) = 1$ \ and
 \ $\EE(Y_0) < \infty$, \ $\EE(|X_0|) < \infty$.
\ Then for all \ $s,t\in\RR_+$ \ with \ $s\leq t$, \ we have
  \begin{align}\label{cond_exp_discrete_Y}
  &\EE(Y_t \mid \cF_s)
     = \ee^{-b(t-s)} Y_s + a \int_s^t \ee^{-b(t-u)} \, \dd u,\\  \label{cond_exp_discrete_X}
  &\EE(X_t \mid \cF_s)
    = X_s + \int_s^t (\alpha - \beta \EE(Y_u \mid \cF_s)) \, \dd u \\\nonumber
  &\phantom{\EE(X_t \mid \cF_s)\,}
     = X_s + \alpha (t - s)
       - \beta Y_s \int_s^t \ee^{-b(u-s)}\,\dd u
       - a \beta \int_s^t \left(\int_s^u \ee^{-b(u-v)} \, \dd v\right) \dd u ,
 \end{align}
 and hence
 \begin{align*}
  \begin{bmatrix}
   \EE(Y_t) \\
   \EE(X_t) \\
  \end{bmatrix}
  = \begin{bmatrix}
     \ee^{-bt} & 0 \\
     - \beta \int_0^t \ee^{-bu} \, \dd u & 1 \\
    \end{bmatrix}
    \begin{bmatrix}
     \EE(Y_0) \\
     \EE(X_0) \\
    \end{bmatrix}
    + \begin{bmatrix}
       \int_0^t \ee^{-bu} \, \dd u & 0 \\
       - \beta \int_0^t \left(\int_0^u \ee^{-bv} \, \dd v \right) \dd u & t \\
      \end{bmatrix}
      \begin{bmatrix}
       a \\
       \alpha \\
      \end{bmatrix} .
 \end{align*}
Consequently, if \ $b \in \RR_{++}$, \ then
 \[
   \lim_{t\to\infty} \EE(Y_t) = \frac{a}{b} , \qquad
   \lim_{t\to\infty} t^{-1} \EE(X_t) = \alpha - \frac{\beta a}{b} ,
 \]
 if \ $b = 0$, \ then
 \[
   \lim_{t\to\infty} t^{-1} \EE(Y_t) = a , \qquad
   \lim_{t\to\infty} t^{-2} \EE(X_t) = - \frac{1}{2} \beta a ,
 \]
 if \ $b \in \RR_{--}$, \ then
 \[
   \lim_{t\to\infty} \ee^{bt} \EE(Y_t) = \EE(Y_0) - \frac{a}{b} , \qquad
   \lim_{t\to\infty} \ee^{bt} \EE(X_t)
   = \frac{\beta}{b} \EE(Y_0) - \frac{\beta a}{b^2} .
 \]
\end{Pro}

Based on the asymptotic behavior of the expectations \ $(\EE(Y_t), \EE(X_t))$
 \ as \ $t \to \infty$, \ we recall a classification of the Heston process
 given by the SDE \eqref{Heston_SDE}, see, Barczy and Pap \cite[Definition 2.3]{BarPap}.

\begin{Def}\label{Def_criticality}
Let \ $(Y_t, X_t)_{t\in\RR_+}$ \ be the unique strong solution of the SDE
 \eqref{Heston_SDE} satisfying \ $\PP(Y_0 \in \RR_+) = 1$.
\ We call \ $(Y_t, X_t)_{t\in\RR_+}$ \ subcritical, critical or supercritical if
 \ $b \in \RR_{++}$, \ $b = 0$ \ or \ $b \in \RR_{--}$, \ respectively.
\end{Def}

In the sequel \ $\stoch$, \ $\distr$ \ and \ $\as$ \ will denote convergence
 in probability, in distribution and almost surely, respectively.

The following result states the existence of a unique stationary distribution
 and the ergodicity for the process \ $(Y_t)_{t\in\RR_+}$ \ given by the first
 equation in \eqref{Heston_SDE} in the subcritical case, see, e.g., Cox et
 al.\ \cite[Equation (20)]{CoxIngRos}, Li and Ma \cite[Theorem 2.6]{LiMa} or
 Theorem 3.1 with \ $\alpha = 2$ \ and Theorem 4.1 in Barczy et al.\
 \cite{BarDorLiPap2}.

\begin{Thm}\label{Ergodicity}
Let \ $a, b, \sigma_1 \in \RR_{++}$.
\ Let \ $(Y_t)_{t\in\RR_+}$ \ be the unique strong solution of the first equation
 of the SDE \eqref{Heston_SDE} satisfying \ $\PP(Y_0 \in \RR_+) = 1$.
\ Then
 \renewcommand{\labelenumi}{{\rm(\roman{enumi})}}
 \begin{enumerate}
  \item
   $Y_t \distr Y_\infty$ \ as \ $t \to \infty$, \ and the distribution of
   \ $Y_\infty$ \ is given by
   \begin{align}\label{Laplace}
    \EE(\ee^{-\lambda Y_\infty})
    = \left(1 + \frac{\sigma_1^2}{2b} \lambda \right)^{-2a/\sigma_1^2} ,
    \qquad \lambda \in \RR_+ ,
   \end{align}
   i.e., \ $Y_\infty$ \ has Gamma distribution with parameters
   \ $2a / \sigma_1^2$ \ and \ $2b / \sigma_1^2$, \ hence
   \[
     \EE(Y_\infty) = \frac{a}{b}, \qquad
     \EE(Y_\infty^2) = \frac{(2a+\sigma_1^2)a}{2b^2} , \qquad
     \EE(Y_\infty^3) = \frac{(2a+\sigma_1^2)(a+\sigma_1^2)a}{2b^3} .
   \]
 \item
  supposing that the random initial value \ $Y_0$ \ has the same distribution
  as \ $Y_\infty$, \ the process \ $(Y_t)_{t\in\RR_+}$ \ is strictly stationary.
 \item
  for all Borel measurable functions \ $f : \RR \to \RR$ \ such that
  \ $\EE(|f(Y_\infty)|) < \infty$, \ we have
  \begin{equation}\label{ergodic}
   \frac{1}{T} \int_0^T f(Y_s) \, \dd s \as \EE(f(Y_\infty)) \qquad
   \text{as \ $T \to \infty$.}
  \end{equation}
 \end{enumerate}
\end{Thm}

In what follows we recall some limit theorems for continuous (local) martingales.
We will use these limit theorems later on for studying the asymptotic
 behaviour of least squares estimators of \ $(a, b, \alpha, \beta)$.
\ First we recall a strong law of large numbers for continuous local
 martingales.

\begin{Thm}{\bf (Liptser and Shiryaev \cite[Lemma 17.4]{LipShiII})}
\label{DDS_stoch_int}
Let \ $\bigl( \Omega, \cF, (\cF_t)_{t\in\RR_+}, \PP \bigr)$ \ be a filtered
 probability space satisfying the usual conditions.
Let \ $(M_t)_{t\in\RR_+}$ \ be a square-integrable continuous local martingale with respect to the
 filtration \ $(\cF_t)_{t\in\RR_+}$ \ such that \ $\PP(M_0 = 0) = 1$.
\ Let \ $(\xi_t)_{t\in\RR_+}$ \ be a progressively measurable process such that
 \ $\PP\big( \int_0^t \xi_u^2 \, \dd \langle M \rangle_u < \infty \big) = 1$, \
 $t \in \RR_+$, \ and
 \begin{align}\label{SEGED_STRONG_CONSISTENCY2}
  \int_0^t \xi_u^2 \, \dd \langle M \rangle_u \as \infty \qquad
  \text{as \ $t \to \infty$,}
 \end{align}
 where \ $(\langle M \rangle_t)_{t\in\RR_+}$ \ denotes the quadratic variation
 process of \ $M$.
\ Then
 \begin{align}\label{SEGED_STOCH_INT_SLLN}
  \frac{\int_0^t \xi_u \, \dd M_u}
       {\int_0^t \xi_u^2 \, \dd \langle M \rangle_u} \as 0 \qquad
  \text{as \ $t \to \infty$.}
 \end{align}
If \ $(M_t)_{t\in\RR_+}$ \ is a standard Wiener process, the progressive
 measurability of \ $(\xi_t)_{t\in\RR_+}$ \ can be relaxed to measurability and
 adaptedness to the filtration \ $(\cF_t)_{t\in\RR_+}$.
\end{Thm}

The next theorem is about the asymptotic behaviour of continuous multivariate
 local martingales, see van Zanten \cite[Theorem 4.1]{Zan}.

\begin{Thm}{\bf (van Zanten \cite[Theorem 4.1]{Zan})}\label{THM_Zanten}
Let \ $\bigl( \Omega, \cF, (\cF_t)_{t\in\RR_+}, \PP \bigr)$ \ be a filtered
 probability space satisfying the usual conditions.
Let \ $(\bM_t)_{t\in\RR_+}$ \ be a $d$-dimensional square-integrable continuous local martingale
 with respect to the filtration \ $(\cF_t)_{t\in\RR_+}$ \ such that
 \ $\PP(\bM_0 = \bzero) = 1$.
\ Suppose that there exists a function \ $\bQ : \RR_+ \to \RR^{d \times d}$
 \ such that \ $\bQ(t)$ \ is an invertible (non-random) matrix for all
 \ $t \in \RR_+$, \ $\lim_{t\to\infty} \|\bQ(t)\| = 0$ \ and
 \[
   \bQ(t) \langle \bM \rangle_t \, \bQ(t)^\top \stoch \bfeta \bfeta^\top
   \qquad \text{as \ $t \to \infty$,}
 \]
 where \ $\bfeta$ \ is a \ $d \times d$ random matrix.
Then, for each $\RR^k$-valued random vector \ $\bv$ \ defined on
 \ $(\Omega, \cF, \PP)$, \ we have
 \[
   (\bQ(t) \bM_t, \bv) \distr (\bfeta \bZ, \bv) \qquad
   \text{as \ $t \to \infty$,}
 \]
 where \ $\bZ$ \ is a \ $d$-dimensional standard normally distributed random
 vector independent of \ $(\bfeta, \bv)$.
\end{Thm}

We note that Theorem \ref{THM_Zanten} remains true if the function \ $\bQ$
 \ is defined only on an interval \ $[t_0, \infty)$ \ with some
 \ $t_0 \in \RR_{++}$.

\section{Existence of LSE based on continuous time observations}
\label{section_LSE_continuous}

First, we define the LSE of \ $(a,b,\alpha,\beta)$ \ based on discrete time observations
 \ $(Y_{\frac{i}{n}},X_{\frac{i}{n}})_{i\in\{0,1,\ldots,\nT\}}$, \ $n\in\NN$, \ $T\in\RR_{++}$ \ (see \eqref{LSE_discrete_sample}) by pointing out that
 the sum appearing in this definition of LSE can be considered as an approximation
 of the corresponding sum of the conditional LSE of \ $(a,b,\alpha,\beta)$ \ based on discrete
 time observations \ $(Y_{\frac{i}{n}},X_{\frac{i}{n}})_{i\in\{0,1,\ldots,\nT\}}$, \ $n\in\NN$, \ \ $T\in\RR_{++}$ \ (which was investigated in Barczy et al. \cite{BarPapSza}).
Then we introduce the LSE of \ $(a,b,\alpha,\beta)$ \ based on continuous time observations
 \ $(X_t)_{t\in[0,T]}$, $T\in\RR_{++}$ \ (see \eqref{LSEab_cont} and \eqref{LSEalphabeta_cont}) as the limit in probability of
 the LSE of \ $(a,b,\alpha,\beta)$ \ based on discrete time observations
 \ $(Y_{\frac{i}{n}},X_{\frac{i}{n}})_{i\in\{0,1,\ldots,\nT\}}$,
 \ $n\in\NN$, \ $T\in\RR_{++}$ \ (see Proposition \ref{LEMMA_LSE_exist}).

A LSE of \ $(a, b, \alpha, \beta)$ \ based on discrete time observations
 \ $(Y_{\frac{i}{n}},X_{\frac{i}{n}})_{i\in\{0,1,\ldots,\nT\}}$, \ $n\in\NN$, \ $T\in\RR_{++}$, \ can be obtained by solving the extremum
 problem
 \begin{align}\label{LSE_discrete_sample}
 \begin{split}
  &\bigl(\ha_{T,n}^{\LSED}, \hb_{T,n}^{\LSED}, \halpha_{T,n}^{\LSED}, \hbeta_{T,n}^{\LSED}\bigr) \\
  &:= \argmin_{(a,b,\alpha,\beta)\in\RR^4}
       \sum_{i=1}^{\nT}
        \left[ \left(Y_{\frac{i}{n}} - Y_{\frac{i-1}{n}} - \frac{1}{n} \left(a - b Y_{\frac{i-1}{n}}\right)\right)^2
              + \left(X_{\frac{i}{n}} - X_{\frac{i-1}{n}} - \frac{1}{n} \left(\alpha - \beta Y_{\frac{i-1}{n}}\right)\right)^2 \right] .
 \end{split}
 \end{align}
Here in the notations the letter \ $\mathrm{D}$ \ refers to discrete time observations.
This definition of LSE can be considered as the corresponding one given in Hu and Long \cite[formula (1.2)]{HuLon2}
 for generalized Ornstein-Uhlenbeck processes driven by \ $\alpha$-stable motions, see also Hu and Long
 \cite[formula (3.1)]{HuLon3}.
For a heuristic motivation of the LSE \eqref{LSE_discrete_sample} based on the discrete observations,
 see, e.g., Hu and Long \cite[page 178]{HuLon1} (formulated for Langevin equations),
 and for a mathematical one, see as follows.
By \eqref{cond_exp_discrete_Y}, for all \ $i\in\NN$, \
 \begin{align*}
  Y_{\frac{i}{n}} - \EE(Y_{\frac{i}{n}} \mid \cF_{\frac{i-1}{n}})
     & = Y_{\frac{i}{n}} - \ee^{-{\frac{b}{n}}} Y_{\frac{i-1}{n}} - a\int_{\frac{i-1}{n}}^{\frac{i}{n}} \ee^{-b({\frac{i}{n}}-u)}\,\dd u
       = Y_{\frac{i}{n}} - \ee^{-{\frac{b}{n}}} Y_{\frac{i-1}{n}} - a\int_0^{\frac{1}{n}} \ee^{-bv}\,\dd v\\
     & = \begin{cases}
            Y_{\frac{i}{n}} - Y_{\frac{i-1}{n}} - \frac{a}{n} & \text{if \ $b=0$,}\\
            Y_{\frac{i}{n}} - \ee^{-{\frac{b}{n}}} Y_{\frac{i-1}{n}} + \frac{a}{b}(\ee^{-{\frac{b}{n}}} - 1) & \text{if \ $b\ne0$.}
         \end{cases}
 \end{align*}
Using first order Taylor approximation of \ $\ee^{-{\frac{b}{n}}}$ \ at \ $b = 0$ \ by \ $1-{\frac{b}{n}}$, \ and that of
 \ $\frac{a}{b}(\ee^{-{\frac{b}{n}}} - 1)$ \ at \ $(a, b) = (0, 0)$ \ by \ $-{\frac{a}{n}}$, \ the random variable
 \ $Y_{\frac{i}{n}} - Y_{\frac{i-1}{n}} - \frac{1}{n} (a - b Y_{\frac{i-1}{n}})$ \
 in the definition \eqref{LSE_discrete_sample} of the LSE of \ $(a,b,\alpha,\beta)$
 can be considered as a first order Taylor approximation of
  \[
   Y_{\frac{i}{n}} - \EE(Y_{\frac{i}{n}} \mid Y_0,X_0,Y_{\frac{1}{n}},X_{\frac{1}{n}},\ldots,Y_{\frac{i-1}{n}},X_{\frac{i-1}{n}}) = Y_{\frac{i}{n}} - \EE(Y_{\frac{i}{n}}\mid \cF_{\frac{i-1}{n}}),
  \]
  which appears in the definition of the conditional LSE of \ $(a,b,\alpha,\beta)$ \ based on discrete time observations
 \ $(Y_{\frac{i}{n}},X_{\frac{i}{n}})_{i\in\{0,1,\ldots,\nT\}}$, \ $n\in\NN$, \ $T\in\RR_{++}$.
\ Similarly, by \eqref{cond_exp_discrete_X}, for all \ $i\in\NN$,
 \begin{align*}
  X_{\frac{i}{n}} - \EE(X_{\frac{i}{n}} \mid\cF_{\frac{i-1}{n}})
     &= X_{\frac{i}{n}} - X_{\frac{i-1}{n}} - {\frac{\alpha}{n}} + \beta Y_{\frac{i-1}{n}} \int_{\frac{i-1}{n}}^{\frac{i}{n}} \ee^{-b(u - \frac{i-1}{n}})\,\dd u
       + a\beta \int_{\frac{i-1}{n}}^{\frac{i}{n}}\left(\int_{\frac{i-1}{n}}^u \ee^{-b(u-v)}\,\dd v\right)\dd u\\
     &= X_{\frac{i}{n}} - X_{\frac{i-1}{n}} - \frac{\alpha}{n} + \beta Y_{\frac{i-1}{n}} \int_0^{\frac{1}{n}} \ee^{-bu}\,\dd u
       + a\beta \int_0^{\frac{1}{n}}\left(\int_0^u \ee^{-bv}\,\dd v\right)\dd u\\
     & = \begin{cases}
            X_{\frac{i}{n}} - X_{\frac{i-1}{n}} - \frac{\alpha}{n} + \frac{\beta}{n} Y_{\frac{i-1}{n}} + \frac{a\beta}{2n^2} & \text{if \ $b=0$,}\\
            X_{\frac{i}{n}} - X_{\frac{i-1}{n}} - \frac{\alpha}{n} + \frac{\beta}{b}(1-\ee^{-\frac{b}{n}}) Y_{\frac{i-1}{n}}
              + \frac{a\beta}{b}\bigl(\frac{1}{n} - \frac{1-\ee^{-\frac{b}{n}}}{b}\bigr)  & \text{if \ $b\ne0$.}
         \end{cases}
 \end{align*}
Using first order Taylor approximation of \ $\frac{a\beta}{2n^2}$ \ at \ $(a, \beta) = (0, 0)$ \ by \ $0$, \ that of
 \ $\frac{\beta}{b}(1-\ee^{-\frac{b}{n}})$ \ at \ $(b, \beta) = (0, 0)$ \ by \ $\frac{\beta}{n}$,
  \ and that of
  \ $\frac{a\beta}{b}\bigl(\frac{1}{n} - \frac{1-\ee^{-\frac{b}{n}}}{b}\bigr)
     = \frac{a\beta}{n^2}\sum_{k=0}^\infty (-1)^k \frac{(b/n)^k}{(k+2)!}$ \ at \ $(a, b, \beta) = (0, 0, 0)$ \ by \ $0$,
 \ the random variable \ $X_{\frac{i}{n}} - X_{\frac{i-1}{n}} - \frac{1}{n} (\alpha - \beta Y_{\frac{i-1}{n}})$ \
  in the definition \eqref{LSE_discrete_sample} of the LSE of \ $(a,b,\alpha,\beta)$ \ can be considered
  as a first order Taylor approximation of
 \[
   X_{\frac{i}{n}} - \EE(X_{\frac{i}{n}} \mid Y_0,X_0,Y_{\frac{1}{n}},X_{\frac{1}{n}},\ldots,Y_{\frac{i-1}{n}},X_{\frac{i-1}{n}} ) = X_{\frac{i}{n}} - \EE(X_{\frac{i}{n}} \mid \cF_{\frac{i-1}{n}}) ,
 \]
 which appears in the definition of the conditional LSE of  \ $(a,b,\alpha,\beta)$ \ based on discrete time observations
 \ $(Y_{\frac{i}{n}},X_{\frac{i}{n}})_{i\in\{0,1,\ldots,\nT\}}$, \ $n\in\NN$, \ $T\in\RR_{++}$.

We note that in Barczy et al. \cite{BarPapSza} we proved  strong consistency and asymptotic normality of conditional LSE of
 \ $(a,b,\alpha,\beta)$ \ based on discrete time observations \ $(Y_i,X_i)_{i\in\{1,\ldots,n\}}$, \ $n\in\NN$, \
 starting the process from some known non-random initial value \ $(y_0,x_0)\in\RR_{++}\times\RR$,
 \ as the sample size \ $n$ \ tends to infinity in the subcritical case.

Solving the extremum problem \eqref{LSE_discrete_sample}, we have
 \begin{align*}
  \bigl(\ha_{T,n}^{\LSED}, \hb_{T,n}^{\LSED}\bigr)
  &= \argmin_{(a,b)\in\RR^2} \sum_{i=1}^{\nT} \left(Y_{\frac{i}{n}} - Y_{\frac{i-1}{n}} - \frac{1}{n}\left(a - b Y_{\frac{i-1}{n}}\right)\right)^2 , \\
  \bigl(\halpha_{T,n}^{\LSED}, \hbeta_{T,n}^{\LSED}\bigr)
  &= \argmin_{(\alpha,\beta)\in\RR^2}
      \sum_{i=1}^{\nT} \left(X_{\frac{i}{n}} - X_{\frac{i-1}{n}} - \frac{1}{n}\left(\alpha - \beta Y_{\frac{i-1}{n}}\right)\right)^2 ,
 \end{align*}
 hence, similarly as on page 675 in Barczy et al. \cite{BarDorLiPap}, we get
 \begin{align}\label{LSEab_discrete}
  \begin{split}
  &\begin{bmatrix}
    \ha_{T,n}^{\LSED} \\[1mm]
    \hb_{T,n}^{\LSED}
   \end{bmatrix}
   = n\begin{bmatrix}
      \nT &   -\sum_{i=1}^{\nT} Y_{\frac{i-1}{n}} \\
      -\sum_{i=1}^{\nT} Y_{\frac{i-1}{n}} & \sum_{i=1}^{\nT} Y_{\frac{i-1}{n}}^2
     \end{bmatrix}^{-1}
     \begin{bmatrix}
      Y_{{\frac{\nT}{n}}} - Y_0\\
      -\sum_{i=1}^{\nT} (Y_{\frac{i}{n}} - Y_{\frac{i-1}{n}}) Y_{\frac{i-1}{n}}
     \end{bmatrix} ,
   \end{split}
  \end{align}
 and
  \begin{align}\label{LSEalphabeta_discrete}
    \begin{split}
  &\begin{bmatrix}
    \halpha_{T,n}^{\LSED} \\[1mm]
    \hbeta_{T,n}^{\LSED}
   \end{bmatrix}
   = n\begin{bmatrix}
      \nT & -\sum_{i=1}^{\nT} Y_{\frac{i-1}{n}} \\
      -\sum_{i=1}^{\nT} Y_{\frac{i-1}{n}} & \sum_{i=1}^{\nT} Y_{\frac{i-1}{n}}^2
     \end{bmatrix}^{-1}
     \begin{bmatrix}
       X_{{\frac{\nT}{n}}} - X_0 \\
      -\sum_{i=1}^{\nT} (X_{\frac{i}{n}} - X_{\frac{i-1}{n}}) Y_{\frac{i-1}{n}}
     \end{bmatrix} ,
     \end{split}
 \end{align}
 provided that the inverse exists, i.e., \ $\nT \sum_{i=1}^{\nT} Y_{\frac{i-1}{n}}^2 > \left(\sum_{i=1}^{\nT} Y_{\frac{i-1}{n}}\right)^2$.
\ By Lemma 3.1 in Barczy et al.~\cite{BarPapSza}, for all \ $n\in\NN$ \ and \ $T\in\RR_{++}$ \ with \ $\nT \geq 2$, \ we have
 \ $\PP\left(\nT \sum_{i=1}^{\nT} Y_{\frac{i-1}{n}}^2 > \left(\sum_{i=1}^{\nT} Y_{\frac{i-1}{n}}\right)^2\right)=1$.

\begin{Pro}\label{LEMMA_LSE_exist}
If \ $a \in \RR_{++}$, \ $b \in \RR$, \ $\alpha,\beta\in\RR$, \ $\sigma_1,\sigma_2\in\RR_{++}$, \ \ $\rho\in(-1,1)$, \ and \ $\PP(Y_0 \in \RR_{++})=1$, \
 then for any \ $T\in\RR_{++}$, \ we have
 \[
  \begin{bmatrix}
    \ha_{T,n}^{\LSED} \\[1mm]
    \hb_{T,n}^{\LSED} \\[1mm]
    \halpha_{T,n}^{\LSED} \\[1mm]
    \hbeta_{T,n}^{\LSED}
  \end{bmatrix}
  \stoch
  \begin{bmatrix}
  \ha_T^{\LSE} \\[1mm]
  \hb_T^{\LSE} \\[1mm]
  \halpha_T^{\LSE} \\[1mm]
  \hbeta_T^{\LSE}
  \end{bmatrix}
  \qquad \text{as \ $n\to\infty$,}
 \]
 where
 \begin{align}\label{LSEab_cont}
  \begin{split}
  \begin{bmatrix}
   \ha_T^{\LSE} \\
   \hb_T^{\LSE}
  \end{bmatrix}
  &:= \begin{bmatrix}
       T  & -\int_0^T Y_s \, \dd s \\
       - \int_0^T Y_s \,\dd s & \int_0^T Y_s^2 \, \dd s
      \end{bmatrix}^{-1}
      \begin{bmatrix}
        Y_T - Y_0 \\
       -\int_0^T Y_s \, \dd Y_s
      \end{bmatrix} \\
  &= \frac{1}{T\int_0^T Y_s^2 \, \dd s - \left(\int_0^T Y_s \,\dd s\right)^2}
     \begin{bmatrix}
        (Y_T - Y_0)\int_0^T Y_s^2 \, \dd s - \int_0^T Y_s \, \dd s \int_0^T Y_s \, \dd Y_s \\
        (Y_T - Y_0)\int_0^T Y_s \, \dd s  - T \int_0^T Y_s \, \dd Y_s
      \end{bmatrix},
  \end{split}
 \end{align}
 and
 \begin{align} \label{LSEalphabeta_cont}
  \begin{split}
  \begin{bmatrix}
   \halpha_T^{\LSE} \\
   \hbeta_T^{\LSE}
  \end{bmatrix}
  &:= \begin{bmatrix}
       T & -\int_0^T Y_s \, \dd s \\
       -\int_0^T Y_s \, \dd s & \int_0^T Y_s^2 \, \dd s
      \end{bmatrix}^{-1}
      \begin{bmatrix}
        X_T - X_0 \\
       - \int_0^T Y_s \, \dd X_s
      \end{bmatrix} \\
  &= \frac{1}{T\int_0^T Y_s^2 \, \dd s - \left(\int_0^T Y_s \,\dd s\right)^2}
     \begin{bmatrix}
        (X_T - X_0)\int_0^T Y_s^2 \, \dd s - \int_0^T Y_s \, \dd s \int_0^T Y_s \, \dd X_s \\
        (X_T - X_0)\int_0^T Y_s \, \dd s  - T \int_0^T Y_s \, \dd X_s
      \end{bmatrix} ,
  \end{split}
 \end{align}
  which exist almost surely, since
 \begin{align}\label{positive_LSE}
  \PP\left( T \int_0^T Y_s^2 \, \dd s
            > \left(\int_0^T Y_s \, \dd s\right)^2 \right)
  = 1 \qquad \text{for all \ $T \in \RR_{++}$.}
 \end{align}
By definition, we call \ $\bigl(\ha_T^{\LSE}, \hb_T^{\LSE}, \halpha_T^{\LSE}, \hbeta_T^{\LSE}\bigr)$ \
 the LSE of \ $(a, b, \alpha, \beta)$ \ based on continuous time observations \ $(X_t)_{t\in[0,T]}$, \ $T\in\RR_{++}$.
\end{Pro}

\noindent{\bf Proof.}
First, we check \eqref{positive_LSE}.
Note that \ $\PP(\int_0^T Y_s \, \dd s < \infty) = 1$ \ and
 \ $\PP(\int_0^T Y_s^2 \, \dd s < \infty) = 1$ \ for all \ $T \in \RR_+$,
 \ since \ $Y$ \ has continuous trajectories almost surely.
For each \ $T \in \RR_{++}$, \ put
 \[
   A_T := \{ \omega \in \Omega
              : \text{$t \mapsto Y_t(\omega)$ \ is continuous and non-negative
                      on \ $[0,T]$} \} .
 \]
Then \ $A_T \in \cF$, \ $\PP(A_T) = 1$, \ and for all \ $\omega \in A_T$,
 \ by the Cauchy--Schwarz's inequality, we have
 \[
   T \int_0^T Y_s(\omega)^2 \, \dd s
   \geq \left(\int_0^T Y_s(\omega) \, \dd s\right)^2 ,
 \]
 and
 \ $T \int_0^T Y_s(\omega)^2 \, \dd s
    - \left(\int_0^T Y_s(\omega) \, \dd s\right)^2
    = 0$
 \ if and only if \ $Y_s(\omega) = K_T(\omega)$ \ for almost
 every \ $s \in [0, T]$ \ with some \ $K_T(\omega) \in \RR_+$. \
\ Hence \ $Y_s(\omega) = Y_0(\omega)$ \ for all \ $s \in [0, T]$ \ if
 \ $\omega \in A_T$ \ and
 \ $T \int_0^T Y_s^2(\omega) \, \dd s - \left(\int_0^T Y_s(\omega) \, \dd s\right)^2 = 0$.
\ Consequently, using that  \ $\PP(A_T)=1$, \ we have
 \begin{align*}
  \PP\left( T \int_0^T Y_s^2 \, \dd s - \left(\int_0^T Y_s \, \dd s\right)^2 = 0 \right)
    &= \PP\left( \left\{T \int_0^T Y_s^2 \, \dd s - \left(\int_0^T Y_s \, \dd s\right)^2 = 0 \right\} \cap A_T\right)\\
    &\leq \PP(Y_s = Y_0,\;\forall \,s\in[0,T]) \leq \PP(Y_T = Y_0) =0,
 \end{align*}
 where the last equality follows by the fact that \ $Y_T$ \ is absolutely continuous
 (see, e.g., Alfonsi \cite[Proposition 1.2.11]{Alf2}) together with the law of total probability.
Hence \ $\PP\Big(T \int_0^T Y_s^2 \, \dd s - \Big(\int_0^T Y_s \, \dd s\Big)^2 = 0 \Big) = 0$, \ yielding
 \eqref{positive_LSE}.

Further, we have
 \[
   \frac{1}{n}
   \begin{bmatrix}
     \nT & -\sum_{i=1}^{\nT} Y_{\frac{i-1}{n}} \\
     -\sum_{i=1}^{\nT} Y_{\frac{i-1}{n}} & \sum_{i=1}^{\nT} Y_{\frac{i-1}{n}}^2
   \end{bmatrix}
   \as \begin{bmatrix}
        T & -\int_0^T Y_s \, \dd s \\
        -\int_0^T Y_s \, \dd s & \int_0^T Y_s^2 \, \dd s
       \end{bmatrix} \qquad \text{as \ $n \to \infty$,}
 \]
 since \ $(Y_t)_{t\in\RR_+}$ \ is almost surely continuous.
By Proposition I.4.44 in Jacod and Shiryaev \cite{JSh} with the Riemann sequence of deterministic subdivisions
 \ $\bigl(\frac{i}{n} \wedge T\bigr)_{i\in\NN}$, \ $n \in \NN$, \ and using the almost sure continuity of
 \ $(Y_t, X_t)_{t\in\RR_+}$, \ we obtain
 \begin{align*}
   \begin{bmatrix}
    Y_{{\frac{\nT}{n}}} - Y_0\\
    -\sum_{i=1}^{\nT} (Y_{\frac{i}{n}} - Y_{\frac{i-1}{n}}) Y_{\frac{i-1}{n}}
   \end{bmatrix}
   &\stoch
   \begin{bmatrix}
    Y_T - Y_0 \\
    -\int_0^T Y_s \, \dd Y_s
   \end{bmatrix} \qquad \text{as \ $n \to \infty$,}  \\
   \begin{bmatrix}
    X_{{\frac{\nT}{n}}} - X_0 \\
    -\sum_{i=1}^{\nT} (X_{\frac{i}{n}} - X_{\frac{i-1}{n}}) Y_{\frac{i-1}{n}}
   \end{bmatrix}
   &\stoch
   \begin{bmatrix}
    X_T - X_0 \\
    -\int_0^T Y_s \, \dd X_s
   \end{bmatrix} \qquad \text{as \ $n \to \infty$.}
 \end{align*}
By Slutsky's lemma, using also \eqref{LSEab_discrete}, \eqref{LSEalphabeta_discrete} and \eqref{positive_LSE}, we obtain the assertion.
\proofend

Note that Proposition \ref{LEMMA_LSE_exist} is valid for all \ $b\in\RR$, \ i.e.,
 not only for subcritical Heston models.

We call the attention that \ $(\ha_T^{\LSE} , \hb_T^{\LSE} , \halpha_T^{\LSE} , \hbeta_T^{\LSE})$ \
 can be considered to be based only on  \ $(X_t)_{t\in[0,T]}$, \ since the process
 \ $(Y_t)_{t\in[0,T]}$ \ can be determined using the observations \ $(X_t)_{t\in[0,T]}$ \
 and the initial value \ $Y_0$, \ see Barczy and Pap \cite[Remark 2.5]{BarPap}.
We also point out that Overbeck and Ryd\'en \cite[formulae (22) and (23)]{OveRyd}
 have already come up with the definition of LSE \ $(\ha_T^{\LSE} , \hb_T^{\LSE} )$ \  of \ $(a,b)$ \
 based on continuous time observations \ $(Y_t)_{t\in[0,T]}$, \ $T \in \RR_{++}$, \ for the CIR process \ $Y$.
\ They investigated only the CIR process \ $Y$, \ so our definitions \eqref{LSEab_cont} and \eqref{LSEalphabeta_cont}
 can be considered as generalizations of formulae (22) and (23) in Overbeck and Ryd\'en \cite{OveRyd}
 for the Heston model \eqref{Heston_SDE}.
Overbeck and Ryd\'en \cite[Theorem 3.4]{OveRyd} also proved that the LSE of
 \ $(a,b)$ \ based on continuous time observations can be approximated in probability
 by conditional LSEs of \ $(a,b)$ \ based on appropriate discrete time observations.

In the next remark we point out that the LSE of \ $(a,b,\alpha,\beta)$ \ given in \eqref{LSEab_cont} and
 \eqref{LSEalphabeta_cont} can be approximated using discrete time observations for \ $X$,
 \ which can be reassuring for practical applications, where data in continuous record is not available.

\begin{Rem}
The stochastic integral \ $\int_0^T Y_s \, \dd Y_s$ \ in \eqref{LSEab_cont} is a measurable function of
 \ $(X_s)_{s\in[0,T]}$ \ and \ $Y_0$.
\ Indeed, for all \ $t\in[0,T]$, \ $Y_t$ \ and \ $\int_0^t Y_s\,\dd s$ \ are measurable functions
 of \ $(X_s)_{s\in[0,T]}$ \ and \ $Y_0$, \ i.e., they can be determined from a sample \ $(X_s)_{s\in[0,T]}$ \
 and \ $Y_0$ \ following from a slight modification of Remark 2.5 in Barczy and Pap \cite{BarPap}
 (replacing \ $y_0$ \ by \ $Y_0$), \ and, by It\^{o}'s formula, we have
 \ $\dd(Y_t^2) = 2 Y_t \,\dd Y_t + \sigma_1^2 Y_t \, \dd t$, \ $t \in \RR_+$, \ implying that
 \ $\int_0^T Y_s \, \dd Y_s = \frac{1}{2}
     \big( Y_T^2 - Y_0^2 - \sigma_1^2 \int_0^T Y_s \, \dd s \big)$, \ $T\in\RR_+$.
\ For the stochastic integral \ $\int_0^T Y_s \, \dd X_s$ \ in \eqref{LSEalphabeta_cont}, we have
 \begin{equation}\label{measurability}
   \sum_{i=1}^\nT Y_{\frac{i-1}{n}} (X_{\frac{i}{n}} - X_{\frac{i-1}{n}})
   \stoch
   \int_0^T Y_s \, \dd X_s \qquad
   \text{as \ $n \to \infty$,}
 \end{equation}
 following from Proposition I.4.44 in Jacod and Shiryaev \cite{JSh} with the
 Riemann sequence of deterministic subdivisions
 \ $\left(\frac{i}{n} \land T\right)_{i\in\NN}$, \ $n \in \NN$.
\ Thus, there exists a measurable function \ $\Phi: C([0,T],\RR)\times\RR\to \RR$
 \ such that \ $\int_0^T Y_s \, \dd X_s = \Phi((X_s)_{s\in[0,T]},Y_0)$, \ since
 the convergence in \eqref{measurability} holds almost surely along a suitable
 subsequence, for each \ $n \in \NN$, \ the members of the sequence in
 \eqref{measurability} are measurable functions of \ $(X_s)_{s\in[0,T]}$ \ and \ $Y_0$,
 \ and one can use Theorems 4.2.2 and 4.2.8 in Dudley \cite{Dud}.
Hence the right hand sides of \eqref{LSEab_cont} and \eqref{LSEalphabeta_cont} are measurable functions of
 \ $(X_s)_{s\in[0,T]}$ \ and \ $Y_0$, \ i.e., they are statistics.
\proofend
\end{Rem}

Using the SDE \eqref{Heston_SDE} and Corollary 3.2.20 in Karatzas and Shreve \cite{KarShr}, one can check that
 \begin{align*}
  \begin{bmatrix}
   \ha_T^{\LSE} - a \\
   \hb_T^{\LSE} - b
  \end{bmatrix}
  &= \begin{bmatrix}
      T & -\int_0^T Y_s \, \dd s \\
      -\int_0^T Y_s \, \dd s  & \int_0^T Y_s^2 \, \dd s
     \end{bmatrix}^{-1}
     \begin{bmatrix}
      \sigma_1 \int_0^T Y_s^{1/2} \, \dd W_s \\
      -\sigma_1 \int_0^T Y_s^{3/2} \, \dd W_s \\
     \end{bmatrix} , \\
  \begin{bmatrix}
   \halpha_T^{\LSE} - \alpha \\
   \hbeta_T^{\LSE} - \beta
  \end{bmatrix}
  &= \begin{bmatrix}
      T & -\int_0^T Y_s \, \dd s \\
      -\int_0^T Y_s \, \dd s  & \int_0^T Y_s^2 \, \dd s
     \end{bmatrix}^{-1}
     \begin{bmatrix}
      \sigma_2 \int_0^T Y_s^{1/2} \, \dd \tW_s\\
      -\sigma_2 \int_0^T Y_s^{3/2} \, \dd \tW_s \\
     \end{bmatrix} ,
 \end{align*}
 provided that \ $T \int_0^T Y_s^2 \, \dd s > \left(\int_0^T Y_s \, \dd s\right)^2$, \ where
 \ $\tW_t:=\varrho W_t + \sqrt{1-\varrho^2} B_t$, \ $t\in\RR_+$, \  and hence
 \begin{align}\label{LSE-}
 \begin{split}
  \ha_T^{\LSE} - a
  &= \frac{\sigma_1 \left(\int_0^T Y_s^{1/2} \, \dd W_s\right)
           \left(\int_0^T Y_s^2 \, \dd s\right)
           - \sigma_1 \left(\int_0^T Y_s \, \dd s\right)
             \left(\int_0^T Y_s^{3/2} \, \dd W_s\right)}
          {T \int_0^T Y_s^2 \, \dd s
           - \left(\int_0^T Y_s \, \dd s\right)^2} , \\
  \hb_T^{\LSE} - b
  &= \frac{\sigma_1 \left(\int_0^T Y_s^{1/2} \, \dd W_s\right)
           \left(\int_0^T Y_s \, \dd s\right)
           - \sigma_1 T \int_0^T Y_s^{3/2} \, \dd W_s}
          {T \int_0^T Y_s^2 \, \dd s
           - \left(\int_0^T Y_s \, \dd s\right)^2} , \\
  \halpha_T^{\LSE} - \alpha
  &= \frac{\sigma_2 \left(\int_0^T Y_s^{1/2} \, \dd \tW_s\right)
           \left(\int_0^T Y_s^2 \, \dd s\right)
           - \sigma_2 \left(\int_0^T Y_s \, \dd s\right)
             \left(\int_0^T Y_s^{3/2} \, \dd \tW_s\right)}
          {T \int_0^T Y_s^2 \, \dd s
           - \left(\int_0^T Y_s \, \dd s\right)^2} , \\
  \hbeta_T^{\LSE} - \beta
  &= \frac{\sigma_2 \left(\int_0^T Y_s^{1/2} \, \dd \tW_s\right)
           \left(\int_0^T Y_s \, \dd s\right)
           - \sigma_2 T \int_0^T Y_s^{3/2} \, \dd \tW_s}
          {T \int_0^T Y_s^2 \, \dd s
           - \left(\int_0^T Y_s \, \dd s\right)^2} ,
 \end{split}
 \end{align}
 provided that
 \ $T \int_0^T Y_s^2 \, \dd s > \left(\int_0^T Y_s \, \dd s\right)^2$.

\section{Consistency and asymptotic normality of LSE}
\label{section_cons_asy}

Our first result is about the consistency of LSE in case of subcritical Heston models.

\begin{Thm}\label{Thm_LSE_cons}
If \ $a, b, \sigma_1, \sigma_2 \in \RR_{++}$, \ $\alpha, \beta \in \RR$,
 \ $\varrho \in (-1, 1)$, \ and \ $\PP((Y_0,X_0) \in \RR_{++}\times\RR)=1$, \ then the LSE of
 \ $(a, b, \alpha, \beta)$ \ is strongly consistent, i.e.,
 \ $\bigl(\ha_T^{\LSE}, \hb_T^{\LSE}, \halpha_T^{\LSE}, \hbeta_T^{\LSE}\bigr)
    \as (a, b, \alpha, \beta)$
 \ as \ $T \to \infty$.
\end{Thm}

\noindent{\bf Proof.}
By Proposition \ref{LEMMA_LSE_exist}, there exists a unique LSE
 \ $\bigl(\ha^{\LSE}_T, \hb^{\LSE}_T, \halpha^{\LSE}_T, \hbeta^{\LSE}_T\bigr)$
 \ of \ $(a, b, \alpha, \beta)$ \ for all \ $T\in\RR_{++}$.
\ By \eqref{LSE-}, we have
 \begin{align*}
  \ha^{\LSE}_T - a
  = \frac{\sigma_1 \cdot \frac{1}{T} \int_0^T Y_s \, \dd s
          \cdot \frac{1}{T} \int_0^T Y_s^2 \, \dd s
          \cdot \frac{\int_0^T Y_s^{1/2} \, \dd W_s}{\int_0^T Y_s \, \dd s}
          - \sigma_1\cdot \frac{1}{T} \int_0^T Y_s \, \dd s
            \cdot \frac{1}{T} \int_0^T Y_s^3 \, \dd s
            \cdot \frac{\int_0^T Y_s^{3/2} \, \dd W_s}{\int_0^T Y_s^3 \, \dd s}}
         { \frac{1}{T} \int_0^T Y_s^2 \, \dd s
              - \left(\frac{1}{T} \int_0^T Y_s \, \dd s\right)^2}
 \end{align*}
 provided that
 \ $\int_0^T Y_s \, \dd s \in\RR_{++}$, \ which holds almost surely, see the proof of Proposition \ref{LEMMA_LSE_exist}.
Since, by part (i) of Theorem \ref{Ergodicity}, \ $\EE(Y_\infty)$,
 $\EE(Y_\infty^2)$, $\EE(Y_\infty^3) \in \RR_{++}$, \ part (iii) of Theorem
 \ref{Ergodicity} yields
 \begin{align*}
  \frac{1}{T} \int_0^T Y_s \, \dd s \as \EE(Y_\infty) , \qquad
  \frac{1}{T} \int_0^T Y_s^2 \, \dd s \as \EE(Y_\infty^2) , \qquad
  \frac{1}{T} \int_0^T Y_s^3 \, \dd s \as \EE(Y_\infty^3)
 \end{align*}
 as \ $T \to \infty$, \ and then
  \begin{align*}
    \int_0^T Y_s \, \dd s \as \infty , \qquad
    \int_0^T Y_s^2 \, \dd s \as \infty , \qquad
    \int_0^T Y_s^3 \, \dd s \as \infty
 \end{align*}
 as \ $T \to \infty$.
\ Hence, by a strong law of large numbers for continuous local
 martingales (see, e.g., Theorem \ref{DDS_stoch_int}), we obtain
 \begin{align*}
  \ha^{\LSE}_T - a
  \as \frac{\sigma_1 \cdot \EE(Y_\infty) \cdot \EE(Y_\infty^2) \cdot 0
            - \sigma_1 \cdot \EE(Y_\infty) \cdot \EE(Y_\infty^3) \cdot 0}
           {\EE(Y_\infty^2) - (\EE(Y_\infty))^2}
      = 0 \qquad \text{as \ $T \to \infty$,}
 \end{align*}
 where for the last step we also used that
 \ $\EE(Y_\infty^2) - (\EE(Y_\infty))^2 = \frac{a\sigma_1^2}{2b^2} \in \RR_{++}$.

Similarly, by \eqref{LSE-},
 \begin{align*}
  \hb^{\LSE}_T - b
  &= \frac{\sigma_1
           \cdot \left(\frac{1}{T} \int_0^T Y_s \, \dd s\right)^2
           \cdot \frac{\int_0^T Y_s^{1/2} \, \dd W_s}{\int_0^T Y_s \, \dd s}
           - \sigma_1
             \cdot \frac{1}{T} \int_0^T Y_s^3 \, \dd s
             \cdot \frac{\int_0^T Y_s^{3/2} \, \dd W_s}{\int_0^T Y_s^3 \, \dd s}}
          {\frac{1}{T} \int_0^T Y_s^2 \, \dd s
           - \left(\frac{1}{T} \int_0^T Y_s \, \dd s\right)^2} \\
  &\as \frac{\sigma_1 \cdot (\EE(Y_\infty))^2 \cdot 0
            - \sigma_1 \cdot \EE(Y_\infty^3) \cdot 0}
           {\EE(Y_\infty^2) - (\EE(Y_\infty))^2}
      = 0 \qquad \text{as \ $T \to \infty$.}
 \end{align*}
One can prove
 \[
   \halpha^{\LSE}_T - \alpha \as 0  \qquad \text{and} \qquad
   \hbeta^{\LSE}_T - \beta \as 0  \qquad \text{as \ $T \to \infty$}
 \]
 in a similar way.
\proofend

Our next result is about the asymptotic normality of LSE in case of subcritical Heston models.

\begin{Thm}\label{Thm_LSE}
If \ $a, b, \sigma_1, \sigma_2 \in \RR_{++}$, \ $\alpha, \beta \in \RR$,
 \ $\varrho \in (-1, 1)$ \ and \ $\PP((Y_0,X_0) \in \RR_{++}\times\RR)=1$, \ then the LSE of
 \ $(a, b, \alpha, \beta)$ \ is asymptotically normal, i.e.,
 \begin{align}\label{LSE_normal}
  T^{\frac{1}{2}}
  \begin{bmatrix}
   \ha_T^{\LSE} - a \\
   \hb_T^{\LSE} - b \\
   \halpha_T^{\LSE} - \alpha \\
   \hbeta_T^{\LSE} - \beta
  \end{bmatrix}
  \distr \cN_4\left(\bzero,
     \bS \otimes
     \begin{bmatrix}
      \frac{(2 a + \sigma_1^2) a}{\sigma_1^2 b}
       & \frac{2 a + \sigma_1^2}{\sigma_1^2} \\
      \frac{2 a + \sigma_1^2}{\sigma_1^2}
       & \frac{2b(a + \sigma_1^2)}{\sigma_1^2 a}
     \end{bmatrix} \right) \qquad
  \text{as \ $T \to \infty$,}
 \end{align}
 where \ $\otimes$ \ denotes the tensor product of matrices, and
 \begin{align*}
  \bS := \begin{bmatrix}
          \sigma_1^2 & \varrho \sigma_1 \sigma_2 \\
          \varrho \sigma_1 \sigma_2 & \sigma_2^2
         \end{bmatrix}.
 \end{align*}
With a random scaling, we have
 \begin{align}\label{LSE_normal_random}
   E_{1,T}^{-\frac{1}{2}}\,
   \bI_2\otimes \begin{bmatrix}
                  (T E_{2,T} - E_{1,T}^2)\big( E_{1,T}E_{3,T} - E_{2,T}^2\big)^{-\frac{1}{2}}  & 0 \\
                  - T & E_{1,T} \\
                \end{bmatrix}
  \begin{bmatrix}
   \ha_T^{\LSE} - a \\
   \hb_T^{\LSE} - b \\
   \halpha_T^{\LSE} - \alpha \\
   \hbeta_T^{\LSE} - \beta
  \end{bmatrix}
  \distr \cN_4\left(\bzero,
     \bS \otimes \bI_2 \right)
 \end{align}
 as \ $T \to \infty$, \ where \ $E_{i,T}:=\int_0^T Y^i_s\,\dd s$, \ $T\in\RR_{++}$, \ $i=1,2,3$.
\end{Thm}

\noindent{\bf Proof.}
By Proposition \ref{LEMMA_LSE_exist}, there exists a unique LSE
 \ $\bigl(\ha^{\LSE}_T, \hb^{\LSE}_T, \halpha^{\LSE}_T, \hbeta^{\LSE}_T\bigr)$ \ of
 \ $(a, b, \alpha, \beta)$.
\ {By \eqref{LSE-},} we have
 \begin{align*}
  \sqrt{T} (\ha^{\LSE}_T - a)
  &= \frac{\frac{1}{T} \int_0^T Y_s^2 \, \dd s \,
           \cdot
           \frac{\sigma_1}{\sqrt{T}} \int_0^T Y_s^{1/2} \, \dd W_s
           - \frac{1}{T} \int_0^T Y_s \, \dd s \,
             \cdot
             \frac{\sigma_1}{\sqrt{T}} \int_0^T Y_s^{3/2} \, \dd W_s}
          {\frac{1}{T} \int_0^T Y_s^2 \, \dd s
           - \left(\frac{1}{T} \int_0^T Y_s \, \dd s\right)^2} , \\
  \sqrt{T} (\hb^{\LSE}_T - b)
  &= \frac{\frac{1}{T} \int_0^T Y_s \, \dd s \,
           \cdot
           \frac{\sigma_1}{\sqrt{T}} \int_0^T Y_s^{1/2} \, \dd W_s
           - \frac{\sigma_1}{\sqrt{T}} \int_0^T Y_s^{3/2} \, \dd W_s}
          {\frac{1}{T} \int_0^T Y_s^2 \, \dd s
           - \left(\frac{1}{T} \int_0^T Y_s \, \dd s\right)^2} , \\
  \sqrt{T} (\halpha^{\LSE}_T - \alpha)
  &= \frac{ \frac{1}{T} \int_0^T Y_s^2 \, \dd s \,
           \cdot
           \frac{\sigma_2}{\sqrt{T}} \int_0^T Y_s^{1/2} \, \dd \tW_s
           - \frac{1}{T} \int_0^T Y_s \, \dd s \,
             \cdot
             \frac{\sigma_2}{\sqrt{T}} \int_0^T Y_s^{3/2} \, \dd \tW_s}
          {\frac{1}{T} \int_0^T Y_s^2 \, \dd s
           - \left(\frac{1}{T} \int_0^T Y_s \, \dd s\right)^2} , \\
  \sqrt{T} (\hbeta^{\LSE}_T - \beta)
  &= \frac{ \frac{1}{T} \int_0^T Y_s \, \dd s \,
           \cdot
           \frac{\sigma_2}{\sqrt{T}} \int_0^T Y_s^{1/2} \, \dd \tW_s
           -  \frac{\sigma_2}{\sqrt{T}} \int_0^T Y_s^{3/2} \, \dd \tW_s}
          {\frac{1}{T} \int_0^T Y_s^2 \, \dd s
           - \left(\frac{1}{T} \int_0^T Y_s \, \dd s\right)^2},
 \end{align*}
 provided that
 \ $T \int_0^T Y_s^2 \, \dd s > \left(\int_0^T Y_s \, \dd s\right)^2$, \ which holds almost surely.
Consequently,
 \begin{align}\nonumber
  \sqrt{T}
  \begin{bmatrix}
   \ha_T^{\LSE} - a \\
   \hb_T^{\LSE} - b \\
   \halpha_T^{\LSE} - \alpha \\
   \hbeta_T^{\LSE} - \beta
  \end{bmatrix}
  & = \frac{1}{\frac{1}{T} \int_0^T Y_s^2 \, \dd s
           - \left(\frac{1}{T} \int_0^T Y_s \, \dd s\right)^2}
       \left(\bI_2\otimes \begin{bmatrix}
                            \frac{1}{T} \int_0^T Y_s^2 \, \dd s & \frac{1}{T} \int_0^T Y_s \, \dd s \\
                            \frac{1}{T} \int_0^T Y_s \, \dd s & 1 \\
                          \end{bmatrix}
       \right)
       \frac{1}{\sqrt{T}}\bM_T\\\label{help1}
  &= \left(\bI_2\otimes \begin{bmatrix}
                            1 & -\frac{1}{T} \int_0^T Y_s \, \dd s \\
                            -\frac{1}{T} \int_0^T Y_s \, \dd s & \frac{1}{T} \int_0^T Y_s^2 \, \dd s \\
                          \end{bmatrix}^{-1}
       \right)
       \frac{1}{\sqrt{T}}\bM_T,
 \end{align}
 provided that
 \ $T \int_0^T Y_s^2 \, \dd s > \left(\int_0^T Y_s \, \dd s\right)^2$, \
 which holds almost surely, where
 \[
   \bM_t := \begin{bmatrix}
             \sigma_1\int_0^t Y_s^{1/2} \, \dd W_s \\[1mm]
             -\sigma_1\int_0^t Y_s^{3/2} \, \dd W_s \\[1mm]
             \sigma_2\int_0^t Y_s^{1/2} \, \dd \tW_s \\[1mm]
             -\sigma_2\int_0^t Y_s^{3/2} \, \dd \tW_s
            \end{bmatrix} , \qquad
   t \in \RR_+ ,
 \]
 is a 4-dimensional square-integrable continuous local martingale due to
 \ $\int_0^t\EE(Y_s)\,\dd s <\infty$ \ and \ $\int_0^t\EE(Y_s^3)\,\dd s <\infty$, $t\in\RR_+$.
\ Next, we show that
 \begin{align}\label{CLT}
  \frac{1}{\sqrt{T}} \bM_T \distr \bfeta \bZ \qquad
  \text{as \ $T \to \infty$,}
 \end{align}
 where \ $\bZ$ \ is a 4-dimensional standard normally distributed random vector and
 \ $\bfeta \in \RR^{4 \times 4}$ \ such that
 \[
   \bfeta \bfeta^\top
   = \bS
     \otimes
     \begin{bmatrix}
      \EE(Y_\infty) & -\EE(Y_\infty^2) \\
      -\EE(Y_\infty^2) & \EE(Y_\infty^3)
     \end{bmatrix} .
 \]
Here the two symmetric matrices on the right hand side are positive definite,
 since \ $\sigma_1,\sigma_2\in\RR_{++}$, \ $\varrho\in(-1,1)$, \ $\EE(Y_\infty)=\frac{a}{b} \in \RR_{++}$ \
 and
 \[
   \EE(Y_\infty) \EE(Y_\infty^3) - (-\EE(Y_\infty^2))^2
      = \frac{a^2\sigma_1^2}{4b^4}(2a+\sigma_1^2) \in \RR_{++},
 \]
 and, so is their Kronecker product.
Hence \ $\bfeta$ \ can be chosen, for instance, as the uniquely defined symmetric
 positive definite square root of the Kronecker product of the two matrices in question.
We have
 \[
   \langle \bM \rangle_t
   = \bS
     \otimes
     \begin{bmatrix}
      \int_0^t Y_s \, \dd s & -\int_0^t Y_s^2 \, \dd s \\
      -\int_0^t Y_s^2 \, \dd s & \int_0^t Y_s^3 \, \dd s
     \end{bmatrix} , \qquad t \in \RR_+ .
 \]
By Theorem \ref{Ergodicity}, we have
 \[
   \bQ(t) \langle \bM \rangle_t \, \bQ(t)^\top
   \as \bS
       \otimes
       \begin{bmatrix}
        \EE(Y_\infty) & -\EE(Y_\infty^2) \\
        -\EE(Y_\infty^2) & \EE(Y_\infty^3)
       \end{bmatrix} \qquad
   \text{as \ $t \to \infty$}
 \]
 with \ $\bQ(t) := t^{-1/2} \bI_4$, \ $t\in\RR_{++}$.
\ Hence, Theorem \ref{THM_Zanten} yields \eqref{CLT}.
Then, by \eqref{help1}, Slutsky's lemma yields
 \begin{align*}
  \sqrt{T}
  \begin{bmatrix}
   \ha_T^{\LSE} - a \\
   \hb_T^{\LSE} - b \\
   \halpha_T^{\LSE} - \alpha \\
   \hbeta_T^{\LSE} - \beta
  \end{bmatrix}
  \distr
  \left(\bI_2\otimes \begin{bmatrix}
                            1 & -\EE(Y_\infty)\\
                            -\EE(Y_\infty) & \EE(Y_\infty^2) \\
                          \end{bmatrix}^{-1}
       \right)
  \bfeta \bZ \distre \cN_4(\bzero,\bSigma)
  \qquad \text{as \ $T \to \infty$,}
 \end{align*}
 where (applying the identities \ $(\bA\otimes\bB)^\top = \bA^\top \otimes \bB^\top$ \ and
  \ $(\bA\otimes\bB)(\bC\otimes\bD) = (\bA\bC)\otimes(\bB\bD)$)
 \begin{multline*}
  \bSigma:=
   \left(\bI_2\otimes
         \begin{bmatrix}
           1 & -\EE(Y_\infty)\\
           -\EE(Y_\infty) & \EE(Y_\infty^2) \\
         \end{bmatrix}^{-1}
       \right)
   \bfeta \EE(\bZ \bZ^\top) \bfeta^\top
   \left(\bI_2\otimes
          \begin{bmatrix}
            1 & -\EE(Y_\infty)\\
            -\EE(Y_\infty) & \EE(Y_\infty^2) \\
           \end{bmatrix}^{-1}
       \right)^\top \\
 \begin{aligned}
  &= \left(\!\bI_2\otimes\!
         \begin{bmatrix}
           1 & -\EE(Y_\infty)\\
           -\EE(Y_\infty) & \EE(Y_\infty^2) \\
         \end{bmatrix}^{-1}
       \right)\!\!
     \left( \!\bS
            \otimes\!
            \begin{bmatrix}
             \EE(Y_\infty) & -\EE(Y_\infty^2) \\
             -\EE(Y_\infty^2) & \EE(Y_\infty^3)
            \end{bmatrix} \right) \!\!
     \left(\!\bI_2\otimes\!
         \begin{bmatrix}
           1 & -\EE(Y_\infty)\\
           -\EE(Y_\infty) & \EE(Y_\infty^2) \\
         \end{bmatrix}^{-1}
       \right) \\
  &= \left( \bI_2
            \bS \bI_2 \right)
     \otimes
     \left( \begin{bmatrix}
             1 & -\EE(Y_\infty)\\
              -\EE(Y_\infty) & \EE(Y_\infty^2) \\
            \end{bmatrix}^{-1}
            \begin{bmatrix}
             \EE(Y_\infty) & -\EE(Y_\infty^2) \\
             -\EE(Y_\infty^2) & \EE(Y_\infty^3)
            \end{bmatrix}
            \begin{bmatrix}
           1 & -\EE(Y_\infty)\\
           -\EE(Y_\infty) & \EE(Y_\infty^2) \\
         \end{bmatrix}^{-1} \right) \\
  &= \frac{1}{( \EE(Y_\infty^2) - (\EE(Y_\infty))^2 )^2} \\
  &\phantom{=\;\;} \times \bS \otimes
     \begin{bmatrix}
      \left(\EE(Y_\infty) \EE(Y_\infty^3) - (\EE(Y_\infty^2))^2\right) \EE(Y_\infty)
       & \EE(Y_\infty) \EE(Y_\infty^3) - (\EE(Y_\infty^2))^2 \\
      \EE(Y_\infty) \EE(Y_\infty^3) - (\EE(Y_\infty^2))^2
       & \EE(Y_\infty^3) - 2 \EE(Y_\infty) \EE(Y_\infty^2) + (\EE(Y_\infty))^3
     \end{bmatrix} ,
 \end{aligned}
 \end{multline*}
 which yields \eqref{LSE_normal}.
Indeed, by Theorem \ref{Ergodicity}, an easy calculation shows that
 \begin{align}\label{help2}
 \begin{split}
  &\left(\EE(Y_\infty) \EE(Y_\infty^3) - (\EE(Y_\infty^2))^2\right) \EE(Y_\infty)
    = \frac{a^3\sigma_1^2}{4b^5}(2a+\sigma_1^2),\\
  &\EE(Y_\infty) \EE(Y_\infty^3) - (\EE(Y_\infty^2))^2
    = \frac{a^2\sigma_1^2}{4b^4}(2a+\sigma_1^2),\\
  &\EE(Y_\infty^3) - 2 \EE(Y_\infty) \EE(Y_\infty^2) + (\EE(Y_\infty))^3
    = \frac{a\sigma_1^2}{2b^3}(a+\sigma_1^2),\\
  &\EE(Y_\infty^2) - \left(\EE(Y_\infty)\right)^2
    = \frac{a\sigma_1^2}{2b^2}.
  \end{split}
 \end{align}

Now we turn to prove \eqref{LSE_normal_random}.
Slutsky's lemma, \eqref{LSE_normal} and \eqref{help2} yield
 \begin{align*}
  & E_{1,T}^{-\frac{1}{2}}\,
   \bI_2\otimes \begin{bmatrix}
                  (T E_{2,T} - E_{1,T}^2)\big( E_{1,T}E_{3,T} - E_{2,T}^2\big)^{-\frac{1}{2}}  & 0 \\
                  - T & E_{1,T} \\
                \end{bmatrix}
   \begin{bmatrix}
   \ha_T^{\LSE} - a \\
   \hb_T^{\LSE} - b \\
   \halpha_T^{\LSE} - \alpha \\
   \hbeta_T^{\LSE} - \beta
   \end{bmatrix} \\
 &\qquad = \oE_{1,T}^{-\frac{1}{2}}\,
    \bI_2\otimes \begin{bmatrix}
                  (\oE_{2,T} - \oE_{1,T}^2)\big( \oE_{1,T}\oE_{3,T} - \oE_{2,T}^2\big)^{-\frac{1}{2}}  & 0 \\
                  - 1 & \oE_{1,T} \\
                \end{bmatrix}
  \sqrt{T}\begin{bmatrix}
      \ha_T^{\LSE} - a \\
      \hb_T^{\LSE} - b \\
      \halpha_T^{\LSE} - \alpha \\
      \hbeta_T^{\LSE} - \beta
   \end{bmatrix}  \\
 &\qquad \distr (\EE(Y_\infty))^{-\frac{1}{2}}
         \bI_2\otimes
              \begin{bmatrix}
                \big(\EE(Y_\infty^2) - (\EE(Y_\infty))^2\big)
                   \big( \EE(Y_\infty)\EE(Y_\infty^3) - (\EE(Y_\infty^2))^2 \big)^{-\frac{1}{2}}  & 0 \\
                 - 1 & \EE(Y_\infty) \\
              \end{bmatrix}\\
  &\phantom{\qquad \distr\;}\times
       \cN_4\left(\bzero,
     \bS \otimes
     \begin{bmatrix}
      \frac{(2 a + \sigma_1^2) a}{\sigma_1^2 b}
       & \frac{2 a + \sigma_1^2}{\sigma_1^2} \\
      \frac{2 a + \sigma_1^2}{\sigma_1^2}
       & \frac{2b(a + \sigma_1^2)}{\sigma_1^2 a}
     \end{bmatrix} \right)\\
   &\phantom{\qquad} \distre  \cN_4(0,\bXi)
       \qquad \text{as \ $T \to \infty$,}
 \end{align*}
 where \ $\oE_{i,T}:=\frac{1}{T}\int_0^T Y_s^i\,\dd s$, \ $T\in\RR_{++}$, \ $i=1,2,3$, \ and,
 applying the identities \ $(\bA\otimes\bB)^\top = \bA^\top \otimes \bB^\top$, \
 \ $(\bA\otimes\bB)(\bC\otimes\bD) = (\bA\bC)\otimes(\bB\bD)$, \ and using \eqref{help2},
 \begin{align*}
  \bXi &:= \frac{1}{\EE(Y_\infty)}
         \left(\bI_2\otimes
              \begin{bmatrix}
                \big(\EE(Y_\infty^2) - (\EE(Y_\infty))^2\big)
                   \big( \EE(Y_\infty)\EE(Y_\infty^3) - (\EE(Y_\infty^2))^2 \big)^{-\frac{1}{2}}  & 0 \\
                 - 1 & \EE(Y_\infty) \\
              \end{bmatrix}\right)\\
      &\phantom{:=\,}
         \times \left(\bS \otimes
         \begin{bmatrix}
           \frac{(2 a + \sigma_1^2) a}{\sigma_1^2 b}
            & \frac{2 a + \sigma_1^2}{\sigma_1^2} \\
              \frac{2 a + \sigma_1^2}{\sigma_1^2}
            & \frac{2b(a + \sigma_1^2)}{\sigma_1^2 a}
         \end{bmatrix}\right)\\
     &\phantom{:=\,}
        \times \left( \bI_2\otimes
             \begin{bmatrix}
                 \big(\EE(Y_\infty^2) - (\EE(Y_\infty))^2\big)
                    \big( \EE(Y_\infty)\EE(Y_\infty^3) - (\EE(Y_\infty^2))^2 \big)^{-\frac{1}{2}}  & 0 \\
                    - 1 & \EE(Y_\infty) \\
               \end{bmatrix} \right)^\top
 \end{align*}
 \begin{align*}
     & = \frac{1}{\EE(Y_\infty)}
          (\bI_2 \bS \bI_2) \otimes
          \Bigg(\begin{bmatrix}
                \big(\EE(Y_\infty^2) - (\EE(Y_\infty))^2\big)
                   \big( \EE(Y_\infty)\EE(Y_\infty^3) - (\EE(Y_\infty^2))^2 \big)^{-\frac{1}{2}}  & 0 \\
                 - 1 & \EE(Y_\infty) \\
         \end{bmatrix}\\
      &\phantom{ = \frac{1}{\EE(Y_\infty)} (\bI_2 \bS \bI_2) \otimes\;}
        \times \begin{bmatrix}
           \frac{(2 a + \sigma_1^2) a}{\sigma_1^2 b}
            & \frac{2 a + \sigma_1^2}{\sigma_1^2} \\
              \frac{2 a + \sigma_1^2}{\sigma_1^2}
            & \frac{2b(a + \sigma_1^2)}{\sigma_1^2 a}
         \end{bmatrix} \\
      &\phantom{ = \frac{1}{\EE(Y_\infty)} (\bI_2 \bS \bI_2) \otimes\;}
         \times \begin{bmatrix}
                  \big(\EE(Y_\infty^2) - (\EE(Y_\infty))^2\big)
                    \big( \EE(Y_\infty)\EE(Y_\infty^3) - (\EE(Y_\infty^2))^2 \big)^{-\frac{1}{2}}  & -1 \\
                   0 & \EE(Y_\infty) \\
              \end{bmatrix} \Bigg)\\
      & =\frac{b}{a} \bS
          \otimes
         \Bigg( \begin{bmatrix}
             \sigma_1(2a+\sigma_1^2)^{-\frac{1}{2}} & 0 \\
            -1 & \frac{a}{b} \\
          \end{bmatrix}
          \begin{bmatrix}
           \frac{(2 a + \sigma_1^2) a}{\sigma_1^2 b}
            & \frac{2 a + \sigma_1^2}{\sigma_1^2} \\
              \frac{2 a + \sigma_1^2}{\sigma_1^2}
            & \frac{2b(a + \sigma_1^2)}{\sigma_1^2 a}
         \end{bmatrix}
         \begin{bmatrix}
             \sigma_1(2a+\sigma_1^2)^{-\frac{1}{2}} & -1 \\
            0 & \frac{a}{b} \\
          \end{bmatrix} \Bigg)\\
       &= \bS\otimes \bI_2.
 \end{align*}
Thus we obtain \eqref{LSE_normal_random}.
\proofend

Next, we formulate a corollary of Theorem \ref{Thm_LSE} presenting separately
 the asymptotic behavior of the LSE of \ $(a,b)$ \ based on continuous time observations
 \ $(Y_t)_{t\in[0,T]}$, $T>0$. \
We call the attention that Overbeck and Ryd\'en \cite[Theorem 3.6]{OveRyd} already derived
 this asymptotic behavior (for more details on the role of the initial distribution, see the Introduction),
 however the covariance matrix of the limit normal distribution in their Theorem 3.6 is somewhat complicated.
It turns out that it can be written in a much simpler form by making a simple reparametrization
 of the SDE (1) in Overbeck and Ryd\'en \cite{OveRyd}, estimating \ $-b$ \ instead of \ $b$ \ (with the notations
 of Overbeck and Ryd\'en \cite{OveRyd}), i.e., considering the SDE \eqref{Heston_SDE} and estimating \ $b$ \ (with our notations).

\begin{Cor}\label{Rem_Ove_Ryd}
If \ $a, b, \sigma_1 \in \RR_{++}$, \ and \ $\PP(Y_0\in \RR_{++})=1$, \ then the LSE of \ $(a, b)$ \
 given in \eqref{LSEab_cont} based on continuous time observations \ $(Y_t)_{t\in[0,T]}$, $T>0$, \ is
 strongly consistent and asymptotically normal, i.e.,
 \ $\bigl(\ha_T^{\LSE}, \hb_T^{\LSE}\bigr)\as (a, b)$ \ as \ $T \to \infty$, \ and
 \begin{align*}
  T^{\frac{1}{2}}
  \begin{bmatrix}
   \ha_T^{\LSE} - a \\
   \hb_T^{\LSE} - b
  \end{bmatrix}
  \distr \cN_2\left(\bzero,
     \begin{bmatrix}
      \frac{(2 a + \sigma_1^2) a}{b}
       & 2 a + \sigma_1^2 \\
      2 a + \sigma_1^2
       & \frac{2b(a + \sigma_1^2)}{a}
     \end{bmatrix} \right) \qquad
  \text{as \ $T \to \infty$.}
 \end{align*}
\end{Cor}

\section{Numerical illustrations}\label{section_numerical}

In this section, first, we demonstrate some methods for the simulation of the Heston model \eqref{Heston_SDE},
 and then we illustrate Theorem \ref{Thm_LSE_cons} and convergence \eqref{LSE_normal} in Theorem \ref{Thm_LSE}
 using generated sample paths of the Heston model \eqref{Heston_SDE}.
We will consider a subcritical Heston model \eqref{Heston_SDE} (i.e., \ $b \in \RR_{++}$) \ with a known
 non-random initial value \ $(y_0, x_0) \in \RR_{++} \times \RR$.
 \ Note that in this case the augmented filtration \ $(\cF_t)_{t\in\RR_+}$ \ corresponding to \ $(W_t,B_t)_{t\in\RR_+}$
 \ and the initial value \ $(y_0,x_0)\in\RR_{++}\times\RR$, \ in fact, does not depend on \ $(y_0,x_0)$.
\ We recall five simulation methods which differ from each other in how the CIR process in the Heston model \eqref{Heston_SDE}
 is simulated.

In what follows, let \ $\eta_{k}$, \ $k \in \{1, \ldots, N\}$, \ be independent standard normally distributed random
 variables with some \ $N \in \NN$, \ and put \ $t_k := k  \frac{T}{N}$, \ $k \in \{0, 1, \ldots, N\}$, \ with some \ $T \in \RR_{++}$.

Higham and Mao \cite{HigMao} introduced the Absolute Value Euler (AVE) method
\[
 Y^{(N)}_{t_{k}}=Y^{(N)}_{t_{k-1}}+(a-bY^{(N)}_{t_{k-1}})(t_{k}-t_{k-1})+\sigma_{1}\sqrt{|Y^{(N)}_{t_{k-1}}|}\sqrt{t_{k}-t_{k-1}}\,\eta_{k} ,
  \qquad k \in \{1, \ldots, N\} ,
\]
 with \ $Y^{(N)}_0=y_0$ \ for the approximation of the CIR process, where \ $a,b,\sigma_1 \in \RR_{++}$.
\ This scheme does not preserve non-negativity of the CIR process.

The Truncated Euler (TE) scheme uses the discretization
 \[
 Y^{(N)}_{t_{k}}=Y^{(N)}_{t_{k-1}}+(a-bY^{(N)}_{t_{k-1}})(t_{k}-t_{k-1})+\sigma_{1}\sqrt{\max(Y^{(N)}_{t_{k-1}},0)}\sqrt{t_{k}-t_{k-1}}\,\eta_{k},
 \qquad k \in \{1, \ldots, N\},
 \]
 with \ $Y^{(N)}_0=y_0$, \ where \ $a,b,\sigma_1 \in\RR_{++}$, \ for approximation of the CIR process \ $Y$, \
 see, e.g., Deelstra and Delbaen \cite{DeDe}.
This scheme does not preserve non-negativity of the CIR process.

The Symmetrized Euler (SE) method gives an approximation of the CIR process \ $Y$ \ via the recursion
 \[
 Y^{(N)}_{t_{k}}=\left|Y^{(N)}_{t_{k-1}}+\left(a-bY^{(N)}_{t_{k-1}}\right)(t_{k}-t_{k-1})
  +\sigma_{1}\sqrt{Y^{(N)}_{t_{k-1}}}\sqrt{t_{k}-t_{k-1}}\,\eta_{k}\right|,
   \qquad k \in \{1, \ldots, N\},
 \]
 with \ $Y^{(N)}_0=y_0$, \ where \ $a,b,\sigma_1 \in\RR_{++}$, \ see, Diop \cite{Dio} or Berkaoui et al. \cite{BeBoDi}
 (where the method is analyzed for more general SDEs including so-called alpha-root processes as well with diffusion coefficient
 \ $\sqrt[\alpha]{x}$ \ with \ $\alpha\in(1,2]$ \ instead of \ $\sqrt{x}$).
\ This scheme gives a non-negative approximation of the CIR process \ $Y$.

The following two methods do not directly simulate the CIR process \ $Y$, \ but its square root \ $Z=(Z_t:=\sqrt{Y_t})_{t\in\RR_+}$.
 \ If \ $a>\frac{\sigma_1^2}{2}$, \ then \ $\PP(Y_t\in\RR_{++}, \;\forall\, t\in\RR_+)=1$, \ and, by It\^{o}'s formula,
  \[
    \dd Z_t = \left( \left( \frac{a}{2} - \frac{\sigma_1^2}{8}\right)\frac{1}{Z_t} - \frac{b}{2}Z_t\right)\dd t
               +\frac{\sigma_1}{2}\,\dd W_t,\qquad t\in\RR_+.
  \]
The Drift Explicit Square Root Euler (DESRE) method (see, e.g., Kloeden and Platen \cite[Section 10.2]{KloPla}
 or Hutzenthaler et al.\ \cite[equation (4)]{HutJenKlo} for general SDEs) simulates
 \ $Z$ \ by
 \[
 Z^{(N)}_{t_k}=Z^{(N)}_{t_{k-1}}+\left(\left(\frac{a}{2}-\frac{\sigma_{1}^2}{8}\right)\frac{1}{Z^{(N)}_{t_{k-1}}}
                 -\frac{b}{2}Z^{(N)}_{t_{k-1}}\right)(t_{k}-t_{k-1})+\frac{\sigma_{1}}{2}\sqrt{t_{k}-t_{k-1}}\,\eta_{k},
                 \qquad k \in \{1, \ldots, N\},
 \]
 with \ $Z^{(N)}_0=\sqrt{y_0}$, \ where  \ $a > \frac{\sigma_{1}^2}{2}$ \ and \ $b,\sigma_1\in\RR_{++}$.
\ Here note that \ $\PP(Z^{(N)}_{t_k} =0)=0$, \ $k\in\{1,\ldots,N\}$, \ since \ $Z^{(N)}_{t_k}$ \ is absolutely continuous.
Transforming back, i.e., \ $Y^{(N)}_{t_k}=(Z^{(N)}_{t_k})^2$, $k\in\{0,1,\ldots,N\}$, \ gives a
 non-negative approximation of the CIR process \ $Y$.

The Drift Implicit Square Root Euler (DISRE) method (see, Alfonsi \cite{Alf} or Dereich et al.\ \cite{DerNeuSzp})
 simulates \ $Z$ \ by
 \[
 Z^{(N)}_{t_k}=Z^{(N)}_{t_{k-1}}+\left(\left(\frac{a}{2}-\frac{\sigma_{1}^2}{8}\right)\frac{1}{Z^{(N)}_{t_k}}
                 -\frac{b}{2}Z^{(N)}_{t_k}\right)(t_{k}-t_{k-1})+\frac{\sigma_{1}}{2}\sqrt{t_{k}-t_{k-1}}\,\eta_{k},
                 \qquad k \in \{1, \ldots, N\},
 \]
 with \ $Z^{(N)}_0=\sqrt{y_0}$, \ where  \ $a > \frac{\sigma_{1}^2}{2}$ \ and \ $b,\sigma_1\in\RR_{++}$.
\ This recursion has a unique positive solution given by
 \begin{align*}
 Z^{(N)}_{t_{k}}=\frac{Z^{(N)}_{t_{k-1}}+\frac{\sigma_{1}}{2}\sqrt{t_{k}-t_{k-1}}\,\eta_{k}}{2+b(t_{k}-t_{k-1})}
          +\sqrt{\frac{\left(Z^{(N)}_{t_{k-1}}+\frac{\sigma_{1}}{2}\sqrt{t_{k}-t_{k-1}}\,\eta_{k}\right)^{2}}{(2+b(t_{k}-t_{k-1}))^{2}}
                                 +\frac{\left(a-\frac{\sigma_{1}^2}{4}\right)(t_{k}-t_{k-1})}{2+ b(t_{k}-t_{k-1})}}
\end{align*}
 for \ $k \in \{1, \ldots, N\}$ \ with \ $Z^{(N)}_0=\sqrt{y_0}$.
\ Transforming again back, i.e., \ $Y^{(N)}_{t_k}=(Z^{(N)}_{t_k})^2$, $k\in\{0,1,\ldots,N\}$, \ gives a strictly
 positive approximation of the CIR process \ $Y$.

We mention that there exist so-called exact simulation methods for the CIR process, see, e.g., Alfonsi \cite[Section 3.1]{Alf2}.
In our simulations, we will use the SE, DESRE and DISRE methods for approximating the CIR process
 which preserve non-negativity of the CIR process.

The second coordinate process \ $X$ \ of the Heston process \eqref{Heston_SDE} will be approximated via the usual
 Euler-Maruyama scheme given by
 \begin{align}\label{EM_scheme_X}
    X^{(N)}_{t_{k}} = X^{(N)}_{t_{k-1}}+(\alpha-\beta Y^{(N)}_{t_{k-1}})(t_{k}-t_{k-1})
               + \sigma_2\sqrt{Y^{(N)}_{t_{k-1}}}\sqrt{t_{k}-t_{k-1}}\big(\varrho\,\eta_{k}+\sqrt{1-\varrho^2}\,\zeta_k\big)
 \end{align}
 for \ $k \in \{1, \ldots, N\}$ \ with \ $X^{(N)}_0=x_0$, \ where \ $\alpha,\beta\in\RR$, \ $\sigma_2\in\RR_{++}$, \ $\varrho\in(-1,1)$,
 \ and \ $\zeta_{k}$, \ $k \in \{1, \ldots, N\}$, \ be independent standard normally distributed random variables independent
 of \ $\eta_k$, \ $k\in\{1,\ldots,N\}$.
\ Note that in \eqref{EM_scheme_X} the factor \ $\sqrt{Y^{(N)}_{t_{k-1}}}$ \ appears, which is well-defined in case of
 the CIR process \ $Y$ \ is approximated by the SE, DESRE or DISRE methods, that we will consider.

We also mention that there exist exact simulation methods for the Heston process \eqref{Heston_SDE}, see, e.g.,
 Broadie and Kaya \cite{BroKay} or Alfonsi \cite[Section 4.2.6]{Alf2}.

We will approximate the estimator \ $\bigl(\ha_T^{\LSE}, \hb_T^{\LSE}, \halpha_T^{\LSE}, \hbeta_T^{\LSE}\bigr)$ \
 given in \eqref{LSEab_cont} and \eqref{LSEalphabeta_cont} using the generated sample paths of \ $(Y,X)$.
\ For this, we need to simulate, for a large time \ $T\in\RR_{++}$, \ the random variables
 \begin{align*}
   Y_T,\quad X_T,\quad I_{1,T}:=\int_0^T Y_s\,\dd s,\quad I_{2,T}:=\int_0^TY_s^2\,\dd s,
   \quad I_{3,T}:=\int_0^T Y_s\,\dd Y_s, \quad I_{4,T}:=\int_0^T Y_s\,\dd X_s.
 \end{align*}
We can easily approximate the \ $I_{i,T}$, $i\in\{1,2,3,4\}$, \ respectively, by
 \begin{align*}
 & I^N_{1,T}:= \sum_{k=1}^N Y^{(N)}_{t_{k-1}} (t_k-t_{k-1})
            = \frac{T}{N} \sum_{k=1}^N Y^{(N)}_{t_{k-1}},
  \qquad I^N_{2,T}:= \sum_{k=1}^N (Y^{(N)}_{t_{k-1}})^2 (t_k-t_{k-1})
            = \frac{T}{N} \sum_{k=1}^N (Y^{(N)}_{t_{k-1}})^2 ,     \\
 & I^N_{3,T}:= \sum_{k=1}^N Y^{(N)}_{t_{k-1}} (Y^{(N)}_{t_k}-Y^{(N)}_{t_{k-1}}),
     \qquad  I^N_{4,T}:= \sum_{k=1}^N Y^{(N)}_{t_{k-1}} (X^{(N)}_{t_k}-X^{(N)}_{t_{k-1}}).
 \end{align*}
Hence, we can approximate \ $\ha_T^{\LSE}$, \ $\hb_T^{\LSE}$, \ $\halpha_T^{\LSE}$, \ and \ $\hbeta_T^{\LSE}$ \ by
 \begin{align*}
    \ha_T^{(N)}:=\frac{(Y^{(N)}_T - y_0) I^N_{2,T} - I^N_{1,T} I^N_{3,T}}{TI^N_{2,T} - (I^N_{1,T})^2},\qquad\qquad
    \hb_T^{(N)}:=\frac{(Y^{(N)}_T - y_0) I^N_{1,T} - T I^N_{3,T}}{TI^N_{2,T} - (I^N_{1,T})^2},\\
    \halpha_T^{(N)}:= \frac{(X^{(N)}_T - x_0) I^N_{2,T} - I^N_{1,T} I^N_{4,T}}{TI^N_{2,T} - (I^N_{1,T})^2},\qquad\qquad
    \hbeta_T^{(N)}:= \frac{(X^{(N)}_T - x_0) I^N_{1,T} - T I^N_{4,T}}{TI^N_{2,T} - (I^N_{1,T})^2}.
 \end{align*}
We point out that \ $\ha_T^{(N)}$, \ $\hb_T^{(N)}$, \ $\halpha_T^{(N)}$ \ and \ $\hbeta_T^{(N)}$ \ are well-defined, since
 \begin{align*}
    TI^N_{2,T} - (I^N_{1,T})^2
      = \frac{T^2}{N} \sum_{k=1}^N \left( Y^{(N)}_{t_k} - \frac{1}{N} \sum_{k=1}^N Y^{(N)}_{t_{k-1}} \right)^2
      \geq 0,
 \end{align*}
 and
 \begin{align*}
    TI^N_{2,T} - (I^N_{1,T})^2 = 0
   & \qquad \Longleftrightarrow\qquad
    Y^{(N)}_{t_k} = \frac{1}{N} \sum_{\ell=1}^N Y^{(N)}_{t_{\ell-1}},
    \qquad k\in\{1,\ldots,N\}\\
   & \qquad \Longleftrightarrow\qquad
    Y^{(N)}_0 = Y^{(N)}_{t_1} = \cdots = Y^{(N)}_{t_{N-1}}.
 \end{align*}
Consequently, using that \ $Y^{(N)}_{t_1}$ \ is absolutely continuous together with the law of total probability, we have
 \ $\PP(TI^N_{2,T} - (I^N_{1,T})^2 \in\RR_{++})=1$.

For the numerical implementation, we take \ $y_0=0.2$, \ $x_0=0.1$, \ $a=0.4$, \ $b=0.3$, \ $\alpha=0.1$, \ $\beta=0.15$,
 \ $\sigma_{1}=0.4$, \ $\sigma_{2}=0.3$, \ $\rho=0.2$, \ $T=3000$, \ and \ $N=30000$ \ (consequently,
 \ $t_{k}-t_{k-1}=0.1$, $k\in\{1,\ldots,N\}$).
\ Note that \ $a>\frac{\sigma_1^2}{2}$ \ with this choice of parameters.
\ We simulate \ $10000$ \ independent trajectories of \ $(Y_T,X_T)$ \ and the normalized error
 \ $T^{\frac{1}{2}}\bigl(\ha_T^{\LSE}-a, \hb_T^{\LSE}-b, \halpha_T^{\LSE}-\alpha, \hbeta_T^{\LSE}-\beta\bigr)$.
\ Table \ref{Table1} contains the empirical mean of \ $Y^{(N)}_T$ \ and \ $\frac{1}{T} X^{(N)}_T$, \ based
 on \ $10000$ \ independent trajectories of \ $(Y_T,X_T)$, \ and the (theoretical) limit
 \ $\lim_{t\to\infty}\EE(Y_t)=\frac{a}{b}$ \ and \ $\lim_{t\to\infty}t^{-1}\EE(X_t)=\alpha-\frac{\beta a}{b}$, \ respectively
 (following from Proposition \ref{Pro_moments}), using the schemes SE, DESRE and DISRE for simulating the CIR process.
\begin{table}[ht!]
\centering
\begin{small}
\begin{tabular}{|c|c|c|c|}
  \hline
   \text{Empirical mean of \ $Y^{(N)}_T$ \ and \ $\frac{1}{T} X^{(N)}_T$} & SE & DESRE & DISRE  \\
  \hline
  $\displaystyle{\lim_{t\to\infty} E(Y_t)=\frac{a}{b} = 1.3333}$ & 1.321025 & 1.325539 & 1.331852 \\
  \hline
  $\displaystyle{\lim_{t\to\infty}t^{-1}E(X_t)=\alpha - \frac{\beta a}{b} = -0.1}$ & -0.09978663 & -0.100054 & -0.09941841 \\
  \hline
\end{tabular}
\end{small}
\caption{Empirical mean of \ $Y^{(N)}_T$ \ (first row) and \ $\frac{1}{T} X^{(N)}_T$ \ (second row).}
\label{Table1}
\end{table}

Henceforth, we will use the above choice of parameters except that \ $T=5000$ \ and \ $N=50000$ \
 (yielding \ $t_{k}-t_{k-1}=0.1$, \ $k\in\{1,\ldots,N\}$).

In Table \ref{Table2} we calculate the expected bias \ ($\EE(\widehat{\theta}^{\mathrm{LSE}}_T-\theta)$),
 \ the \ $L_{1}$-norm of error \ ($\EE|\widehat{\theta}^{\mathrm{LSE}}_T-\theta|$) \ and the \ $L_{2}$-norm of error
 \ \big($\big(\EE(\widehat{\theta}^{\mathrm{LSE}}_T-\theta)^{2}\big)^{1/2}$\big), \ where \ $\theta\in\{a,b,\alpha,\beta\}$,
 \ using the scheme DISRE for simulating the CIR process.
\begin{table}[ht!]
\centering
\begin{small}
\begin{tabular}{|c|c|c|c|}
  \hline
    Errors & Expected bias & $L_1$-norm of error & $L_2$-norm of error \\
  \hline
  $a$ & -0.01089369 & 0.0153848 & 0.0190123 \\
  \hline
  $b$ & -0.007639168 & 0.01189344 & 0.01474495 \\
  \hline
  $\alpha$ & 0.0001779072 & 0.00957648 & 0.0120646 \\
  \hline
  $\beta$ & 0.0001402452 & 0.007776999 & 0.009771835 \\
  \hline
\end{tabular}
\end{small}
\caption{Expected bias, \ $L_1$- and \ $L_2$-norm of error using DISRE scheme.}
\label{Table2}
\end{table}

In Table \ref{Table3} we give the relative errors \ $(\widehat{\theta}^{(N)}_T - \theta)/\theta$, \ where
 \ $\theta\in\{a,b,\alpha,\beta\}$, \ for \ $T=5000$ \ using the scheme DISRE for simulating the CIR process.
 \begin{table}[ht!]
\centering
\begin{small}
\begin{tabular}{|c|c|}
  \hline
  Relative errors & $T=5000$   \\
  \hline
  $(\widehat{a}^{(N)}_T- a)/a$ & -0.02723421 \\
  \hline
  $(\widehat{b}^{(N)}_T- b)/b$ & -0.02546389 \\
  \hline
  $(\widehat{\alpha}^{(N)}_T - \alpha)/\alpha$ & 0.001779072 \\
  \hline
  $(\widehat{\beta}^{(N)}_T - \beta)/\beta$ & 0.0009349683 \\
  \hline
\end{tabular}
\end{small}
\caption{Relative errors using DISRE scheme.}
\label{Table3}
\end{table}

\begin{figure}[ht!]
 \centering
 \includegraphics[width=16cm,height=9cm]{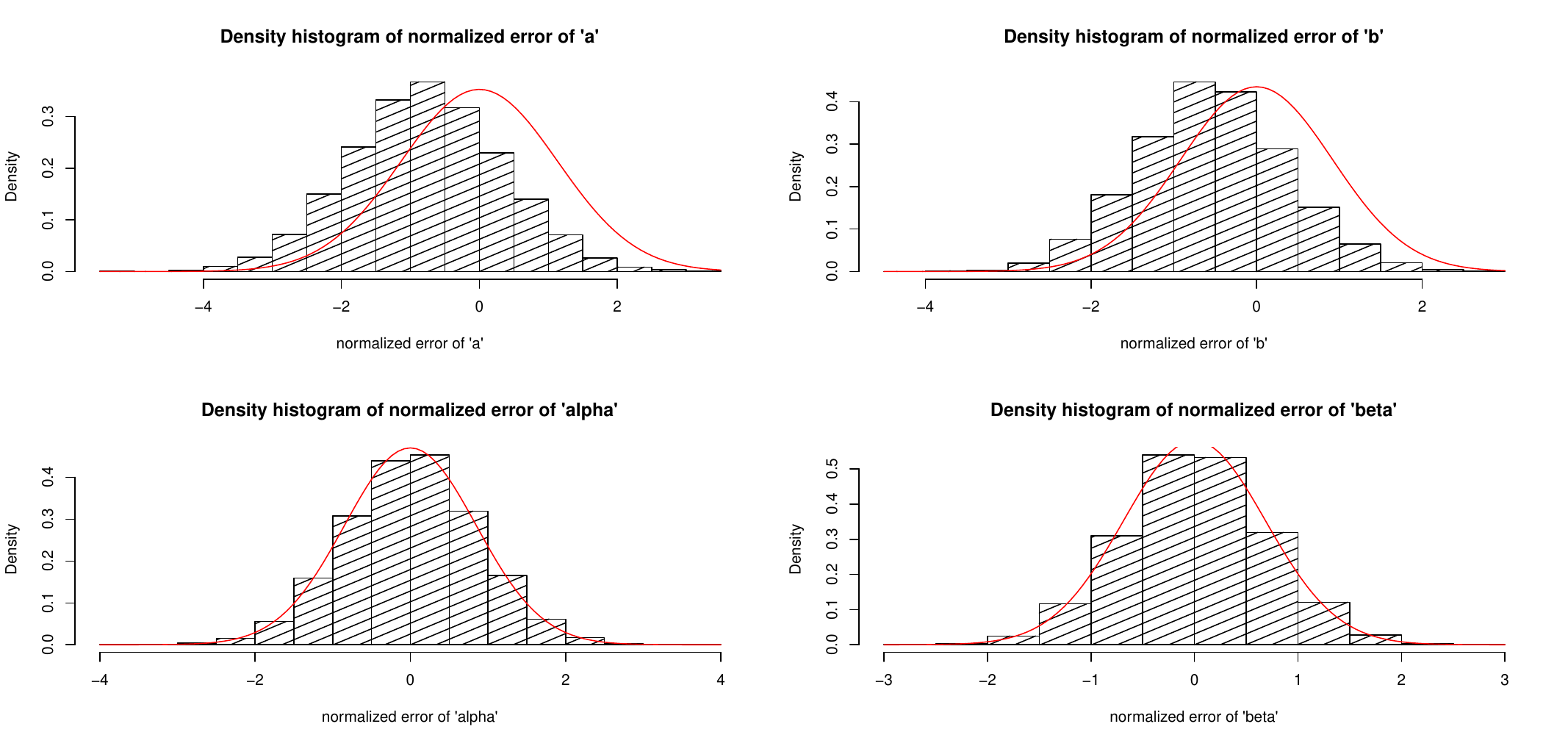}
 \caption{ In the first line from left to right, the density histograms of the
  normalized errors of $T^{1/2}(\ha_T^{(N)} - a)$ and \ $T^{1/2}(\hb_T^{(N)} - b)$, in the second line from left to right, the density histograms of the normalized errors of \ $T^{1/2}(\halpha_T^{(N)} - \alpha)$ \
   and $T^{1/2}(\hbeta_T^{(N)} - \beta)$. \ In each case, the red line denotes the density function of the
   corresponding normal limit distribution. }
   \label{fig}
 \end{figure}

In Figure \ref{fig}, we illustrate the limit law of each coordinate of the LSE
 \ $\bigl(\ha_T^{\LSE}, \hb_T^{\LSE}, \halpha_T^{\LSE}, \hbeta_T^{\LSE}\bigr)$ \ given in \eqref{LSE_normal}.
To do so, we plot the obtained density histograms of each of its coordinates based on \ $10000$ \ independently
 generated trajectories using the scheme DISRE for simulating the CIR process, we also plotted
 the density functions of the corresponding normal limit distributions in red.
With the above choice of parameters, as a consequence of \eqref{LSE_normal}, we have
 \begin{align*}
  &T^{\frac{1}{2}}(\ha_T^{\LSE} - a) \distr \cN\left(0,\frac{a}{b}(2a+\sigma_1^2)\right) = \cN(0,1.28)
     \qquad \text{as \ $T\to\infty$,}\\
  &T^{\frac{1}{2}}(\hb_T^{\LSE} - b) \distr \cN\left(0,\frac{2b}{a}(a+\sigma_1^2)\right) = \cN(0,0.84)
     \qquad \text{as \ $T\to\infty$,}\\
  &T^{\frac{1}{2}}(\halpha_T^{\LSE} - \alpha) \distr \cN\left(0,\frac{a\sigma_2^2}{b\sigma_1^2}(2a+\sigma_1^2)\right) = \cN(0,0.72)
     \qquad \text{as \ $T\to\infty$,}\\
  &T^{\frac{1}{2}}(\hbeta_T^{\LSE} - \beta) \distr \cN\left(0,\frac{2b\sigma_2^2}{a\sigma_1^2}(a+\sigma_1^2)\right) = \cN(0,0.4725)
     \qquad \text{as \ $T\to\infty$.}
 \end{align*}
In case of the parameters \ $a$ \ and \ $b$, \ one can see a bias in Figure \ref{fig}, which, in our
 opinion, may be related with the different speeds of weak convergence for the LSE of \ $(a,b)$ \ and that
 of \ $(\alpha,\beta)$, \ and with the bad performance of the applied discretization scheme for \ $Y$.

Table \ref{Table4} contains the skewness and excess kurtosis of \ $T^{\frac{1}{2}}(\widehat{\theta}^{(N)}_T - \theta)$, \ where
 \ $\theta\in\{a,b,\alpha,\beta\}$, \ using the scheme DISRE for simulating the CIR process.
This confirms our results in \eqref{LSE_normal} as well.
 \begin{table}[ht!]
\centering
\begin{small}
\begin{tabular}{|c|c|c|c|c|}
  \hline
   Skewness and excess kurtosis& $T^{\frac{1}{2}}(\ha^{(N)}_{T}-a)$ & $T^{\frac{1}{2}}(\hb^{(N)}_{T}-b)$ & $T^{\frac{1}{2}}(\halpha^{(N)}_{T}-\alpha)$ & $T^{\frac{1}{2}}(\hbeta^{(N)}_{T}-\beta)$ \\
  \hline
  Skewness & 0.04915124 & 0.04544189 & -0.02317407 & -0.01399869 \\
  \hline
  Excess kurtosis & 0.07666643 & 0.05226811 & 0.09994108 & 0.07877347 \\
  \hline
\end{tabular}
\end{small}
\caption{Skewness and excess kurtosis using the scheme DISRE for simulating the CIR process.}
\label{Table4}
\end{table}

Using the Anderson-Darling and Jarque-Bera tests, we test whether
 each of the coordinates of \ $T^{\frac{1}{2}}\bigl(\ha_T^{\LSE} - a, \hb_T^{\LSE} - b, \halpha_T^{\LSE} - \alpha,
  \hbeta_T^{\LSE} - \beta\bigr)$ \ follows a normal distribution or not for \ $T=5000$.
\ In Table \ref{Table5} we give the test values and (in paranthesis) the p-values of the Anderson-Darling and Jarque-Bera tests
 using the scheme DISRE for simulating the CIR process
 (the \ $*$ \ after a p-value denotes that the p-value in question is greater than any
 reasonable signifance level).
It turns out that, with this choice of parameters, \
at any reasonable significance level the Anderson-Darling test accepts that \ $T^{\frac{1}{2}}(\ha_T^{\LSE} - a)$,
 \ $T^{\frac{1}{2}}(\hb_T^{\LSE} - b)$, \ $T^{\frac{1}{2}}(\halpha_T^{\LSE} - \alpha)$, \
 and \ $T^{\frac{1}{2}}(\hbeta_T^{\LSE} - \beta)$ \ follow normal laws.
The Jarque-Bera test also accepts that \ $T^{\frac{1}{2}}(\hb_T^{\LSE} - b)$, \ $T^{\frac{1}{2}}(\halpha_T^{\LSE} - \alpha)$,
 \ and \ $T^{\frac{1}{2}}(\hbeta_T^{\LSE} - \beta)$ \ follow normal laws,
 but rejects that \ $T^{\frac{1}{2}}(\ha_T^{\LSE} - a)$ \ follows a normal law.
 \begin{table}[ht!]
\centering
\begin{small}
\begin{tabular}{|c|c|c|c|c|}
  \hline
  Test of normality & $T^{\frac{1}{2}}(\ha^{(N)}_{T}-a)$ & $T^{\frac{1}{2}}(\hb^{(N)}_{T}-b)$ & $T^{\frac{1}{2}}(\halpha^{(N)}_{T}-\alpha)$ & $T^{\frac{1}{2}}(\hbeta^{(N)}_{T}-\beta)$ \\
  \hline
  Anderson-Darling & 0.34486 ($0.4857^{*}$) & 0.62481 ($0.1037^{*}$) & 0.34078 ($0.4962^{*}$) & 0.35232 ($0.467^{*}$) \\
  \hline
  Jarque-–Bera & 6.5162 (0.03846) & 4.6077 ($0.09987^{*}$) & 5.1089 ($0.07774^{*}$) & 2.9528 ($0.2285^{*}$) \\
  \hline
\end{tabular}
\end{small}
\caption{Test of normalilty in case of \ $y_0=0.2$, \ $x_0=0.1$, \ $a=0.4$, \ $b=0.3$, \ $\alpha=0.1$, \ $\beta=0.15$,
 \ $\sigma_{1}=0.4$, \ $\sigma_{2}=0.3$, \ $\rho=0.2$, \ $T=5000$, \ and \ $N=50000$ \ generating \ $10000$ \ independent sample paths
 using the scheme DISRE for simulating the CIR process.}
\label{Table5}
\end{table}

All in all, our numerical illustrations are more or less in accordance with our theoretical results in
 \eqref{LSE_normal}.
Finally, we note that we used the open source software \ R \ for making the simulations.

\end{document}